\documentclass{amsproc}

\usepackage{html}		

\newtheorem{theorem}{Theorem}[section]
\newtheorem{lemma}[theorem]{Lemma}
\newtheorem{proposition}[theorem]{Proposition}
\newtheorem{corollary}[theorem]{Corollary}

\theoremstyle{definition}
\newtheorem{definition}[theorem]{Definition}

\theoremstyle{remark}
\newtheorem{remark}[theorem]{Remark}

\numberwithin{equation}{section}

\newcommand{\abs}[1]{\lvert#1\rvert}

\newcommand{\C}{\ensuremath{\mathbb{C}}}
\newcommand{\Q}{\ensuremath{\mathbb{Q}}}
\newcommand{\R}{\ensuremath{\mathbb{R}}}
\newcommand{\Z}{\ensuremath{\mathbb{Z}}}
\newcommand{\N}{\ensuremath{\mathbb{N}}}

\newcommand{\Ccut}{\ensuremath{\mathbb{C}^\prime}}
\newcommand{\HH}{\ensuremath{\mathbb{H}}} 

\newcommand{\RR}{\ensuremath{\mathcal{R}}} 
\newcommand{\JJ}{\ensuremath{\mathcal{J}}} 
\newcommand{\PP}{\ensuremath{\mathrm{P}}}
\newcommand{\tTT}{\ensuremath{\tilde{\mathrm{T}}}} 
 
\newcommand{\hH}{\ensuremath{\mathrm{H}}} 
\newcommand{\thH}{\ensuremath{\tilde{\mathrm{H}}}}
 
\newcommand{\FE}{\ensuremath{\mathrm{FE}}} 
\newcommand{\Gmod}{\ensuremath{{\Gamma}}}
\newcommand{\nmid}{\ensuremath{{\not|}\;}}

\newcommand{\re}[1]{\ensuremath{{\mathrm{Re}\!\left( #1 \right)}}}
\newcommand{\im}[1]{\ensuremath{{\mathrm{Im}\!\left( #1 \right)}}}

\newcommand{\SL}[1]{\ensuremath{{\mathrm{SL}\!\left(2, #1 \right)}}}
\newcommand{\GL}[1]{\ensuremath{{\mathrm{GL}\!\left(2, #1 \right)}}}

\newcommand{\Gnull}[1]{\ensuremath{{\Gamma_0\left(#1\right)}}}
\newcommand{\Matrix}[4]{{\textstyle  \left( {#1\atop #3} \: {#2 \atop  #4} \right)}}

\begin{document}

\title[Hecke algebra]{A realization of the Hecke algebra  on the space of period functions for $\Gamma_0(n)$}

\author{M.~Fraczek}
\address{Institut f\"ur Theoretische Physik, Technische Univerist\"at Clausthal, Clausthal-Zellerfeld, Germany}
\email{marekf@gmx.net}
\author{D.~Mayer}
\address{Institut f\"ur Theoretische Physik, Technische Univerist\"at Clausthal, Clausthal-Zellerfeld, Germany}
\email{dieter.mayer@tu-clausthal.de}
\author{T.~M\"uhlenbruch}
\address{Institut f\"ur Theoretische Physik, Technische Univerist\"at Clausthal, Clausthal-Zellerfeld, Germany}
\email{tobias.muehlenbruch@tu-clausthal.de}
\thanks{The second and third author were supported by the Deutsche Forschungsgemeinschaft through the DFG Research Project "Transfer operators and non arithmetic quantum chaos" (Ma~633/16-1)}
\thanks{\today}

\date{\today}



\begin{abstract}
The standard realization of the Hecke algebra on classical holomorphic cusp forms and the corresponding period polynomials is well known.
In this article we consider a nonstandard realization of the Hecke algebra on Maass cusp forms for the Hecke congruence subgroups $\Gnull{n}$.
We show that the vector valued period functions derived recently by Hilgert, Mayer and Movasati as special eigenfunctions of the transfer operator for $\Gnull{n}$ are indeed related to the Maass cusp forms for these groups.
This leads also to a simple interpretation of the ``Hecke like'' operators of these authors in terms of the aforementioned nonstandard realization of the Hecke algebra on the space of vector valued period functions.
\end{abstract}

\maketitle


\section{Introduction}
\label{A}

There are basically two approaches to the theory of period functions attached to Maass cusp forms for Fuchsian groups: one is just an extension of the Eichler, Manin, Shimura theory of period polynomials for holomorphic cusp forms.
Thereby the Maass cusp forms are related to the period functions by a certain integral transformation as discussed for $\SL{\Z}$ in \cite{LZ01}, and for $\Gnull{n}$ in \cite{Mu05}.

Another, in a certain sense dynamical, approach starts from the geodesic flow on the corresponding surface of constant negative curvature and its transfer operator.
It identifies the period polynomials and period functions as certain eigenfunctions of the analytically continued transfer operator of this flow \cite{HMM05}.
This second approach obviously is also of some interest in the theory of quantum chaos. Indeed, eigenstates of a quantum system, namely a particle moving freely on a surface of constant negative curvature with the hyperbolic Laplace-Beltrami operator as its Schroedinger operator, are thereby related to classical objects, namely eigenfunctions of the classical transfer operator.
Such an exact connection of quantum states with functions attached to the classical system cannot be established up to now within the more familiar approach to quantum chaos through the Selberg-Gutzwiller trace formula \cite{Sa95}.

For arithmetic Fuchsian groups like the modular group and its subgroups there exists a whole family of symmetries for the quantum system described by the so called Hecke operators.
They commute with each other and the Laplacian, and their existence leads to interesting statistical properties of the spectra of such systems \cite{Sa95}.

In a recent paper the authors of \cite{HMM05} constructed for the Hecke congruence subgroups $\Gnull{n}$ certain linear operators $\tTT_{n,m}$ on the space of eigenfunctions of the corresponding transfer operator which they called``Hecke like'' operators.
The operators were derived by using only the structure of the transfer operators for the groups $\Gnull{nm}$ respectively the closely related Lewis functional equations obeyed by their eigenfunctions.
It is known that in the case of the full modular group $\SL{\Z}$ the "Hecke like" operators $\tTT_{1,p}$ for $p$ prime coincide with the ordinary Hecke operators $H_p$ when acting on the period functions (see \cite{Mu05}, \cite{MM05}).

In the present paper we complete and extend this result to arbitrary groups $\Gnull{n}$ by showing that the operators $\tTT_{n,m}$ of \cite{HMM05} indeed define a certain realization of the Hecke algebra on the space of eigenfunctions of the transfer operator for these groups.

This result follows from a direct relation between the Maass cusp forms and the special eigenfunctions of the transfer operators for the subgroups $\Gnull{nm}$ of $\Gnull{n}$ which the authors in \cite{HMM05} used in their derivation of the Hecke like operators.
Another ingredient in our proof is a certain non standard realization of the Hecke algebra on the Maass cusp forms for $\Gnull{n}$ based on a result by Atkin and Lehner derived in \cite{AL70} in the context of holomorphic cusp forms.

In detail this paper is organized as follows: after recalling in Chapter~\ref{section 1} briefly the construction of the Hecke like operators of \cite{HMM05} via the eigenfunctions of the transfer operator we discuss in Chapter~\ref{B} several cosets of Hecke congruence subgroups and their subgroups.
In Chapter~\ref{periodfunctions} we introduce vector valued Maass cusp forms and the integral transformation leading to the vector valued period functions for $\Gnull{n}$.
We show how the eigenfunctions  of the transfer operators constructed in \cite{HMM05} can be interpreted as integral transforms of certain old Maass cusp forms.
In Chapter~\ref{section3} we introduce a realization of the Hecke algebra on Maass cusp forms which differs slightly from the regular realization in the literature.
Transferring this realization via the aforementioned integral transformation to the space of vector valued period function we see that the Hecke like operators of \cite{HMM05} indeed coincide with this realization.
In Chapter~\ref{section6} we prove our main Theorem stated in Chapter~\ref{section 1}.

\section{The ``Hecke like'' operators of Hilgert, Mayer and Movasati}
\label{section 1}

Let us fix some notations which we will use throughout this paper.
We denote by $\HH=\{x+iy:\, y>0\}$ the hyperbolic plane with the hyperbolic metric induced by $ds^2= \frac{dx^2+dy^2}{y^2}$.
We call the group $\SL{\Z}$ the full modular group.
A group $\Gamma \subset \SL{\Z}$ of finite index $\mu=[\SL{\Z}:\Gamma]$ is called a modular group.
A particular class of modular groups are the \emph{Hecke congruence subgroups} $\Gnull{n}$, $n \in \N$, defined as
\begin{equation}
\Gnull{n} = \left\{ \Matrix{a}{b}{c}{d} \in \SL{\Z}; \, c \equiv 0 \bmod n \right\}.
\end{equation}
Obviously, we have $\Gnull{1}= \SL{\Z}$.
The surfaces $M_\Gamma = \Gamma \backslash \HH$ with $\Gamma$ a modular group are called modular surfaces.
They are covering surfaces of $M=\SL{\Z} \backslash \HH$.
We denote by $T$ respectively $S$  the generators
\begin{equation}
\label{PQ}
T=\Matrix{1}{1}{0}{1}, \qquad S=\Matrix{0}{-1}{1}{0}
\end{equation}
of $\SL{\Z}$.
Lateron we shall also need the matrices
\begin{equation}
\label{M}
M=\Matrix{0}{1}{1}{0}, \qquad T^\prime = MTM =\Matrix {1}{0}{1}{1}.
\end{equation}

Let us introduce the \emph{vector valued period functions} for the modular group $\Gamma$.
These are vector valued functions $\vec{\phi}=(\phi_i)_{1 \leq i \leq \mu}$ with the following properties:
\begin{itemize}
\item 
Each component $\phi_i:\Ccut \to \C$ of $\vec{\phi}$ is holomorphic in the complex cut-plane $\Ccut := \C \setminus (-\infty,0]$.
\item
$\vec{\phi}$ fulfills the three-term functional equation
\begin{equation}
\label{three-term 1}
\vec{\phi}(z) - \chi_{\overline{\Gamma}}(T^{-1})\, \vec{\phi}(z+1) - (z+1)^{-2\beta} \, \chi_{\overline{\Gamma}}\big({T^\prime}^{-1}\big)\, \vec{\phi}\left(\frac{z}{z+1}\right) = 0,
\end{equation}
the so called \emph{generalized Lewis equation}.
Here $\chi_{\overline{\Gamma}}$ is the representation of $\GL{\Z}$ by $\mu \times \mu$ permutation matrices induced from the trivial representation of the subgroup $\overline{\Gamma}$ which is the extension of $\Gamma$ to $\GL{\Z}$ by adjoining the element $\Matrix{1}{0}{0}{-1}$ and, if not yet contained in $\Gamma$, also $\Matrix{-1}{0}{0}{-1}$. We assume thereby that $\Gamma \Matrix{1}{0}{0}{-1} = \Matrix{1}{0}{0}{-1}\Gamma. $
\item
The components $\phi_i$ of $\vec{\phi}$ fulfill certain growth properties for $z \to 0$ and $z \to \infty$ depending on $\beta$ as discussed in \cite{Mu05} and \cite{LZ01}.
\end{itemize}
For $\Gamma$ the Hecke congruence subgroup $\Gnull{n}$ we denote its space of vector valued period functions $\vec{\phi}$ with spectral parameter $\beta$ by $\FE(n,\beta)$.

\begin{remark}
For the full modular group $\SL{\Z}$ the period functions fulfill the original Lewis equation
\begin{equation}
\label{Lewis equation}
\phi(z) - \phi(z+1) - (z+1)^{-2\beta} \, \phi\left(\frac{z}{z+1}\right) = 0
\end{equation}
introduced in \cite{LZ01}.
\end{remark}

Consider the ring $\RR = \Z[\mathrm{Mat}_\ast(2,\Z)]$ of finite linear combinations of $2 \times 2$ integer matrices with nonzero determinant, and the right $\RR$-ideal $\JJ$ with $\JJ := (\mathbf{1} - T - T^\prime)\RR$ and $\mathbf{1}$ the unit matrix.
Moreover, denote by $\RR^+$ the subset $\RR^+ = \Z[\mathrm{Mat}_\ast^+(2,\Z)]\subset \RR$ of finite linear combinations of $2 \times 2$ integer matrices with nonnegative entries and put $\JJ^+ =  \JJ \cap \RR^+$.

We use the familiar \emph{slash action $\big|_\beta$} on functions $f:\Ccut \to \C$ with
\begin{equation}
\label{slash action}
\left(f \big|_\beta h \right)(z) := |ad-bc|^\beta\,(cz+d)^{-2\beta} \, f(hz)
\end{equation}
for $hz= \frac{az+b}{cz+d}$ and $h=\Matrix{a}{b}{c}{d} \in \mathrm{Mat}^+_\ast(2,\Z)$.
In \cite{HMM05} it is shown that this slash action is indeed well defined for complex $\beta$; for $\beta \in 2\Z$ the slash action is defined for all $h \in \mathrm{Mat}_\ast(2,\Z)$.
Obviously the action in (\ref{slash action}) extends  linearly to an action of $\RR^+$ and it can be generalized to an action on vector valued functions $\vec{f} =(f_i)_{1 \leq i \leq \mu}$ with $f_i:\Ccut \to \C$ through
\begin{equation}
\label{slash action 2}
\vec{f} \big|_\beta  h:=  \left( f_i\big|_\beta h \right)_{1 \leq i \leq \mu}
\end{equation}
for elements $h \in \RR^+$.
The Lewis Equation (\ref{three-term 1}) can then be written in the  form
\begin{equation}
\label{three-term 2}
\vec{\phi} - \chi_{\overline{\Gamma}}(T^{-1})\, \vec{\phi}\big|_\beta T - \chi_{\overline{\Gamma}}\big({T^\prime}^{-1}\big)\, \vec{\phi}\big|_\beta T^\prime = \vec{0}.
\end{equation}
We do not know how to solve this equation in general, but it is possible to describe special solutions.

Let $\vec{\psi} = (\psi_i)_{1 \leq i \leq \mu}$ be a vector of elements $\psi_i$ in $\RR^+$ solving the
vector valued equation:
\begin{equation}
\label{three-term 3}
\vec{\psi}- \chi_{\overline{\Gamma}}(T^{-1})\, \vec{\psi} T - \chi_{\overline{\Gamma}}\big({T^\prime}^{-1}\big)\, \vec{\psi} T^\prime =\vec{0}\mod\JJ^+,
\end{equation}
where $\vec{\psi} h:= (\psi_i h)_{1 \leq i \leq \mu}$.
Then one has the obvious
\begin{lemma}
Given any solution $\phi=\phi(z)$ of the Lewis Equation (\ref{Lewis equation}) for $\SL{\Z}$ the functions
\[
\phi_i=\phi_i(z):= \phi\big|_\beta \psi_i (z),
\quad 
i \in \{1, \ldots,\mu\}
\]
solve the generalized Lewis Equation (\ref{three-term 1}) for the modular group $\Gamma$ if the elements $\psi_i \in \RR^+$ solve Equation (\ref{three-term 3}).
\end{lemma}

There is a straightforward solution of Equation (\ref{three-term 3}) given by $\psi_i =\mathbf{1}$, $i \in \{1,\ldots,\mu\}$.
This leads to the special but trivial solution $\phi_i(z)=\phi(z)$, $i \in \{1,\ldots,\mu\}$, of Equation (\ref{three-term 1}).
That this solution exists is not surprising since we know in the case of the full modular group from the work of Lewis and Zagier \cite{LZ01} and for general modular groups from the work of Deitmar and Hilgert \cite{DH04} that there is a 1-1 correspondence between their period functions and their Maass cusp forms  for $\re{\beta} > 0$.
But each Maass cusp form for $\SL{\Z}$ is trivially a Maass cusp form for any modular group $\Gamma$ and hence each period function for $\SL{\Z}$ should also give rise to such a function for any of its modular subgroups.
Consistently with the terminology for automorphic forms we  call the above solution of the Lewis equation for the modular group $\Gamma$ an ``old solution''.

In the following we will restrict our discussion of nontrivial solutions of equations (\ref{three-term 1}), respectively (\ref{three-term 3}), to the Hecke congruence subgroups $\Gnull{n}$ as presented in \cite{HMM05}.
For this we need a special characterization of the index set $\overline{\Gnull{n}} \backslash \GL{\Z}$ as given there.
Consider on $\Z \times \Z$ the equivalence relation $\sim_n$ defined as
\begin{eqnarray*}
(x,y) &\sim_n& (x^\prime, y^\prime) 
\quad \mbox{iff} \quad 
\exists k \in \Z, \gcd(k,n)=1 \mbox{ such that } \\
&& \qquad \qquad \qquad \qquad x^\prime \equiv kx \mod n, y^\prime \equiv ky \mod n
\end{eqnarray*}
together with the natural right action of $\GL{\Z}$
\begin{equation}
\label{right action}
(x,y)\Matrix{a}{b}{c}{d} = (xa+yc, xb+yd),
\end{equation}
obviously compatible with $\sim_n$.
Hence $\GL{\Z}$ acts also on $(\Z \times \Z)_n:= (\Z \times \Z) / \sim_n$.
Denote the elements of $(\Z \times \Z)_n$ by $[x:y]_n$.
It is easy to see that the stabilizer in $\GL{\Z}$ of the element $[0:1]_n \in (\Z \times \Z)_n$ is just the subgroup $\overline{\Gnull{n}}$.
Therefore the following map $\overline{\pi}_n: \overline{\Gnull{n}} \backslash \GL{\Z} \to (\Z \times \Z)_n$ is well defined and injective:
\begin{equation}
\label{pi map}
\overline{\pi}_n \left( \overline{\Gnull{n}} \Matrix{a}{b}{c}{d} \right) 
:= 
[0:1]_n \Matrix{a}{b}{c}{d} = [c:d]_n.
\end{equation}
Denote by $I_n$ with 
\begin{equation}
\label{set In 1}
I_n := \overline{\pi}_n \left(  \overline{\Gnull{n}} \backslash \GL{\Z} \right).
\end{equation}
the image of $\overline{\pi}_n$
It is not very difficult to show \cite{HMM05} that
\begin{equation}
\label{set In 2}
I_n = \left\{ [x:y]_n \in (\Z \times \Z)_n: \, \gcd(x,y,n)=1 \right\}.
\end{equation}

Consider next the subset $P_n \subset \Z \times \Z$ with
\begin{equation}
\label{set Pn 1}
P_n= \left\{ (c,b) \in \Z \times \Z: \, c \geq 1, c |n , 0 \leq b \leq \frac{n}{c}-1, \gcd(c,b,\frac{n}{c})=1 \right\}.
\end{equation}
There is a bijection between $I_n$ and $P_n$ given by the map \cite{HMM05}
\begin{equation}
\label{map Pn In}
P_n \ni (c,b) \longmapsto [c:d_n(c,b)]_n \in I_n
\end{equation}
with 
\begin{equation}
\label{map dn}
d_n(c,b) =  \min_{0 \leq k \leq c-1} \left\{ c+b+k\frac{n}{c}: \gcd(c,b+k\frac{n}{c}) =1 \right\}.
\end{equation}
For simplicity we denote in the following the elements of $I_n$ also by $i$.
The bijection in (\ref{map Pn In}) allows us to identify each element $i \in I_n$ uniquely with a matrix $A_i \in \mathrm{Mat}_n(2,\Z)$ with
\begin{equation}
\label{matrix Ai}
A_i = \Matrix{c}{b}{0}{\frac{n}{c}}
\end{equation}
where $\mathrm{Mat}_n(2,\Z)$ denotes the $2 \times 2$ matrices with integer entries and determinant $n$.
In the following we need certain sets of $2 \times 2$ matrices with nonnegative integer entries:
\begin{eqnarray}
\label{set Sn}
&& S_n = \left\{ \Matrix{a}{b}{c}{d} \in \mathrm{Mat}_n(2,\Z): \, a > c \geq 0, d > b \geq 0 \right\}, \\
\label{set Xn}
&& X_n = \left\{ \Matrix{c}{b}{0}{\frac{n}{c}} \in S_n\right\}
,\qquad 
Y_n = \left\{ \Matrix{\frac{n}{c}}{0}{b}{c} \in S_n\right\}, \mbox{ and}
\\
\label{set Xn star}
&& X_n^\star = \left\{ \Matrix{c}{b}{0}{\frac{n}{c}} \in X_n: \, \gcd(c,b,\frac{n}{c})=1 \right\}.
\end{eqnarray}
Obviously the matrix $A_i$ in (\ref{matrix Ai}) belongs to $X_n^\star$ for all $i \in I_n$.

Next consider the map
\begin{equation}
\label{map K 1}
K: S_n \setminus Y_n \longrightarrow S_n \setminus X_n
\end{equation}
defined as
\begin{equation}
\label{map K 2}
K\left(\Matrix{a}{b}{c}{d} \right) := 
\Matrix{-c+ \lceil \frac{d}{b} \rceil a}{-d+ \lceil \frac{d}{b} \rceil b}{a}{b} 
\end{equation}
with $\lceil r \rceil \in \Z$ determined for $r \in \R$ by $\lceil r \rceil -1 < r \leq \lceil r \rceil$.

There exists for every $A \in S_n \setminus Y_n$ an integer $k_A  > 0$ such that $K^jA \notin Y_n$ for $0 \leq j < k_A$ and $K^{k_A}A \in Y_n$.
For $A \in Y_n$ put $k_A=0$ so that $k_A$ is well defined for all $A \in S_n$.

The following Theorem has been proven in \cite{HMM05}.
\begin{theorem}
\label{theorem 1}
The matrices $\psi_i:= \sum_{j=0}^{k_i} K^j(A_i)$, $i \in I_n$, with $k_i=k_{A_i}$ and $A_i$ as in (\ref{matrix Ai}), solve the Lewis Equation (\ref{three-term 3}) for the group $\Gnull{n}$.
\end{theorem}

As an immediate Corollary one gets
\begin{corollary}
For $\phi=\phi(z)$ any solution of the Lewis Equation (\ref{Lewis equation}) for $\SL{\Z}$ the function $\vec{\phi} =(\phi_i)_{i \in I_n}$ with $\phi_i(z) := \phi\big|_\beta \psi_i$, $i \in I_n$ solves the Lewis Equation (\ref{three-term 1}) for the group $\Gnull{n}$ with the same parameter $\beta$.
\end{corollary}

This result can be generalized in the following way.
Since 
\[
\SL{\Z} = \Gnull{1} \supset \Gnull{n} \supset \Gnull{nm}
\]
for fixed $n$ and all $m = 1,2,\ldots$ one has a natural projection
\begin{equation}
\label{projection}
\sigma_{m,n}: I_{nm} \longrightarrow I_n
\end{equation}
induced from the map $\Gnull{nm}\Matrix{a}{b}{c}{d} \mapsto \Gnull{n} \Matrix{a}{b}{c}{d}$.

In \cite{HMM05} one finds
\begin{lemma}
\label{lemma 1}
For any $i \in I_{nm}$ and $0 \leq j \leq k_{\sigma_{n,m}(i)}$ there exists a unique index $l_{i,j} \in I_n$ such that 
\[
A_{l_{i,j}} \big( K^j (A_{\sigma_{n,m}(i)}) \big) A_i^{-1} \in \SL{\Z}.
\]
\end{lemma}

This on the other hand allowed the authors in \cite{HMM05} to prove
\begin{theorem}
\label{theorem 2}
For $\vec{\psi} =\big( \psi_i\big)_{i\in I_n}$ respectively $\vec{\phi} = \big(\phi_i(z)\big)_{i \in I_n}$ any solution of the Lewis equations (\ref{three-term 3}) respectively (\ref{three-term 1}) with parameter $\beta$ for the group $\Gnull{n}$ the matrices $\vec{\Psi} = \big(\Psi_j\big)_{j\in I_{nm}}$ respectively the functions $\vec{\Phi} = \big(\Phi_j(z)\big)_{j\in I_{nm}}$ with
\begin{equation}
\label{matrices Psi}
\Psi_j:= \sum_{s=0}^{k_{\sigma_{n,m}(j)}} \psi_{l_{j,s}} \, K^s (A_{\sigma_{n,m}(j)}),
\end{equation}
respectively
\begin{equation}
\label{functions Phi}
\Phi_j (z):= \sum_{s=0}^{k_{\sigma_{n,m}(j)}} \Big(\phi_{l_{j,s}} \big|_\beta K^s (A_{\sigma_{n,m}(j)}) \Big) (z),
\end{equation}
solve the corresponding Lewis equations for the group $\Gnull{nm}$ with the same parameter $\beta$.
\end{theorem}

An immediate Corollary is
\begin{corollary}
For $\vec{\phi} = \big( \phi_i(z)\big)_{i\in I_n}$ a vector valued period function for $\Gnull{n}$ the function $\vec{\Phi} = \big(\Phi_j(z)\big)_{j \in I_{nm}}$ in (\ref{functions Phi}) is a vector valued period function for $\Gnull{nm}$ with the same parameter $\beta$.
\end{corollary}

\begin{remark}
The function $\vec{F} = \big( F_j(z)\big)_{j\in I_{nm}}$ with $F_j(z) = \Phi_j(z+1)$ is an eigenfunction of the transfer operator for the group $\Gnull{nm})$ with eigenvalue $\lambda= \pm 1$ if the function $\vec{f}=\big(f_i(z)\big)_{i \in I_n}$ with $f_i(z) = \phi_i(z+1)$ is an eigenfunction of the operator for the group $\Gnull{n}$ with the same eigenvalue $\lambda=\pm 1$.
\end{remark}

Another result in \cite{HMM05} which we need later on is
\begin{proposition}
\label{prop 1}
If the function $\vec{\Phi} = \big( \Phi_i(z) \big)_{i \in I_{nm}}$ solves the Lewis Equation (\ref{three-term 1}) for $\Gnull{nm}$ then the function $\vec{\phi} = \big( \phi_j(z) \big)_{j \in I_n}$ with $\phi_j(z) := \sum_{i \in \sigma_{m,n}^{-1}(j)} \Phi_i(z)$ solves this equation for $\Gnull{n}$.
If, on the other hand, $\vec{\phi} = \big(\phi_j(z)\big)_{j \in I_n}$ solves Equation (\ref{three-term 1}) for $\Gnull{n}$, then $\vec{\Phi} = \big( \Phi_i(z)\big)_{i \in I_{nm}}$ with $\Phi_i(z):= \phi_{\sigma_{m,n}(i)}(z)$ solves this equation for $\Gnull{nm}$.
\end{proposition}

\begin{remark}
The second part of this Proposition shows that any period function for $\Gnull{n}$ determines a ``trivial'' old period function for $\Gnull{nm}$ whose components are just given by the components of the former one.
\end{remark}

An immediate consequence in the case $n=1$, that is for the full modular group $\SL{\Z}$, is 
\begin{corollary}
If $\vec{\Psi} = \big(\Psi_i\big)_{i \in I_m}$ solves Equation (\ref{three-term 3}) for $\Gnull{m}$ then the element $\psi^{(m)}$ with
\begin{equation}
\label{matrix psi(m) 1}
\psi^{(m)} := \sum_{i \in I_m} \Psi_i
\end{equation}
solves this equation for the group $\SL{\Z}$ for the same parameter $\beta$.
\end{corollary}

A straightforward calculation \cite{HMM05} shows that for any $m$ prime one has indeed
\begin{equation}
\label{matrix psi(m) 2}
\psi^{(m)} = \sum_{A \in S_m} A.
\end{equation}

\begin{corollary}
For any period function $\phi= \phi(z)$ for $\SL{\Z}$ and for any $m \in \N$ the function $\tilde{\phi}= \tilde{\phi}(z) := \phi\big|_\beta \psi^{(m)} (z)$ is again a period function for this group with the same parameter $\beta$.
\end{corollary}

In complete analogy one derives from Theorem~\ref{theorem 2} and Proposition~\ref{prop 1} 

\begin{theorem}
\label{theorem 3}
For any vector valued period function $\vec{\phi} = \big( \phi_i \big)_{i \in I_n}$ for $\Gnull{n}$ and any $m \in \N$ the function $\vec{\tilde{\phi}} = \tTT_{n,m} \vec{\phi}$ with
\begin{equation}
\label{Hecke operator 1}
\left( \tTT_{n,m} \vec{\phi} \right)_i (z) 
=
\sum_{s \in \sigma_{m,n}^{-1}(i)} \sum_{j=0}^{k_{\sigma_{n,m}(s)}} \Big(\phi_{l_{s,j}} \big|_\beta K^j(A_{\sigma_{n,m}(s)})\Big) (z)
\end{equation}
is again a period function for $\Gnull{n}$ with the same parameter $\beta$.
\end{theorem}
In particular $\tTT_{n,1}$ is the trivial map $\vec{\phi} \mapsto \vec{\phi}$.

Hence, by using only properties of the transfer operators for the geodesic flows on modular surfaces, the authors in \cite{HMM05} constructed linear operators $\tTT_{n,m}$ mapping the space of vector valued period functions for $\Gnull{n}$ with parameter $\beta$ into itself.
In the case $n=1$ and $m$ prime the operator $\tTT_{1,m}$ reduces to the form
\begin{equation}
\label{Hecke operator 2}
\left( \tTT_{1,m} \phi \right) (z) =\Big( \phi \big|_\beta  \sum_{A \in S_m}  A \Big) (z)
\end{equation}
and hence coincides exactly with the $m^\mathrm{th}$ Hecke operator $\thH_m$ in the form derived by M\"uhlenbruch in \cite{Mu04} for period functions of Maass cusp forms for $\SL{\Z}$.
For $m,n \in N$, $m$ prime and $\gcd(n,m)=1$, the relation between $\tTT_{n,m}$ and Hecke operators given in the form as in \cite{AL70} is discussed in \cite{MM05}.

In the present paper we relate the operators $\tTT_{n,m}$ on the space of period functions to some Hecke operators on the space of Maass cusp forms for arbitrary $n,m \in \N$.
This  allows us to prove the following Theorem:
\begin{theorem}
\label{main theorem}
For fixed $n \in \N$ the operators $\tTT_{n,m}$, $m \in \N$, defined in Theorem~\ref{theorem 3} satisfy
\begin{eqnarray*}
\tTT_{n,p} \, \tTT_{n,p^e}
&=&
\left\{\begin{array}{ll}
\tTT_{n,p^{e+1}}
	& \mbox{for prime } p \mid n, e \in \N\\
\tTT_{n,p^{e+1}} + p \Matrix{p}{0}{0}{p} \tTT_{n,p^{e-1}}
	& \mbox{for prime } p \nmid n, e>1,\\
\tTT_{n,p^2} + (p+1) \Matrix{p}{0}{0}{p} \tTT_{n,1} \quad
	& \mbox{for prime } p \nmid n, e=1 \mbox{ and},
\end{array}\right.\\
\tTT_{n,m} \, \tTT_{n,m^\prime}
&=&
\tTT_{n,mm^\prime}
\qquad\qquad\qquad\qquad\qquad\quad\; \mbox{for } (m,m^\prime)=1.
\end{eqnarray*}
\end{theorem}

Theorem~\ref{main theorem} shows in particular that for any $n$ the family of operators $\{\tTT_{n,m}\}$ is a realization of the Hecke algebra on vector valued period functions.
This realization is slightly different from the standard one as given for instance in \cite{Mi89}.

\section{Cosets of Hecke congruence subgroups and their representatives}
\label{B}

We did not succeed to prove Theorem \ref{main theorem} directly from the definition of the operators $\tilde{T}_{n,m}$ in (\ref{Hecke operator 1}).
Instead we are going to relate in a first step the solutions $\vec{\Phi}=(\Phi_j(z))_{j\in I_{nm}}$ in Theorem \ref{theorem 2} to Maass cusp forms for the group $\Gnull{nm}$ respectively the solutions $\vec{\phi}=(\phi_i(z))_{i\in I_n}$ in Proposition \ref{prop 1} to Maass cusp forms for the group $\Gnull{n}$. This allows us in a second step to relate the operators $\tilde{T}_{n,m}$ to certain Hecke operators on these cusp forms fulfilling commutation relations similar to the ones in Theorem \ref{main theorem}.
For this we need some properties of the  representatives of different cosets of the Hecke congruence subgroups.
\begin{lemma}
\label{B1}
For $B_{m^e}=\Matrix{m^e}{0}{0}{1}$ and $\Gamma_0(n,m)$ the subgroup
\begin{equation}
\label{B1a}
\Gamma_0(n,m) := \left\{ \Matrix{a}{b}{c}{d} \in \SL{\Z} : \,c \equiv 0 \bmod n, b \equiv 0 \bmod m  \right\}
\end{equation}
one has:
\begin{itemize}
\item $\Gamma_0(n,1) = \Gnull{n}$,
\item $\Gamma_0(n,m^{e+l}) B_{m^e} = B_{m^e} \Gamma_0(nm^e,m^l)$ for all $e,l\in \N_0$ and in particular
\item $\Gamma_0(n,m) B_m = B_m \Gnull{nm}$.
\end{itemize}
\end{lemma}

\begin{proof}
By definition $\Gamma_0(n,1) = \Gnull{n}$.

We show the inclusion ``$\subset$'' for the second equality:
take $\Matrix{a}{b}{c}{d} \in \Gamma_0(n,m^{e+l})$,
then
\[
\Matrix{a}{b}{c}{d} \, B_{m^e} = \Matrix{a}{b}{c}{d} \, \Matrix{m^e}{0}{0}{1} = \Matrix{m^ea}{b}{m^ec}{d} = \Matrix{m^e}{0}{0}{1}\, \Matrix{a}{\frac{b}{m^e}}{m^ec}{d} = B_{m^e}\, \Matrix{a}{\frac{b}{m^e}}{m^ec}{d}.
\]
Since $m^{e+l} \mid b$ and $n \mid c$ we find  $\Matrix{a}{\frac{b}{m^e}}{m^ec}{d} \in \Gamma_0(nm^e,m^l)$.

Next consider the inclusion ``$\supset$'': take $\Matrix{a}{b}{c}{d} \in \Gamma_0(nm^e,m^l)$,
then
\[
B_{m^e}\, \Matrix{a}{b}{c}{d}  = \Matrix{m^ea}{m^eb}{c}{d} = \Matrix{a}{m^eb}{\frac{c}{m^e}}{d} \, B_{m^e}.
\]
Since $m^l \mid b$ and $nm^e \mid c$ we find  $\Matrix{a}{m^eb}{\frac{c}{m^e}}{d} \in \Gamma_0(n,m^{e+l})$.
\end{proof}

Denote by $I_{nm,n}$ the index set $I_{nm,n}=\{1,\ldots,[\Gnull{n}:\Gnull{nm}]\}$.
Obviously, the sets $I_{n,1}$ and $I_n$ can be identified.

Next we show the following 
\begin{lemma}
\label{B7}
For $n,m \in \N$ with $\gcd(n,m)=1$ and $g \in \Gnull{n}$ the relation $B_m \, g \, B_m^{-1} \in \Gnull{n,m}$ implies $g \in \Gnull{nm}$.
\end{lemma}

\begin{proof}
Write $g = \Matrix{a}{b}{nc}{d}$.
We have
\[
B_m \, g \, B_m^{-1}
=
\Matrix{m}{0}{0}{1}\Matrix{a}{b}{nc}{d}\Matrix{\frac{1}{m}}{0}{0}{1}
=
\Matrix{a}{mb}{\frac{nc}{m}}{d} \in \Gamma_0(n,m).
\]
Since $\gcd(n,m)=1$ we have $m | c$ and in particular $nm \mid nc$.
Hence $g \in \Gnull{nm}$ 
\end{proof}

\smallskip

For $n,m\in \mathbb{N}$ let $R_j^{nm,n}$, $j \in I_{nm,n}$, be a system of representatives of the right cosets in $\Gnull{nm}\setminus \Gnull{n}$ satisfying
\begin{equation}
\label{2 system of representatives}
\bigsqcup_{j  \in I_{nm,n}} \Gnull{nm}\, R_j^{nm,n} = \Gnull{n}.
\end{equation}
Here $\bigsqcup_{j  \in I_{nm,n}}$ denotes the disjoint union:
\begin{eqnarray*}
&&
\bigcup_{j  \in I_{nm,n}} \Gnull{nm}\, R_j^{nm,n} = \Gnull{n}
\qquad \mbox{and}\\
&&
\Gnull{nm}\, R_{j_1}^{nm,n} \cap \Gnull{nm}\, R_{j_2}^{nm,n} = \emptyset
\qquad \mbox{for all distinct } j_1,j_2 \in I_{nm,n}.
\end{eqnarray*}

\begin{lemma}
\label{B8}
\begin{itemize}
\item
For arbitrary $p,m,n \in \N$ consider a system of representatives $R_j$, $j \in J$ a suitable index set, of the cosets in $\Gamma_0(nm,p) \setminus \Gamma_0(n,p)$.
For distinct $j_1, j_2 \in J$ we have 
\[
\Gamma_0(nm) \, R_{j_1} \cap \Gamma_0(nm) \, R_{j_2} = \emptyset.
\]
\item
For $p$ prime, $p \mid n$, let $R_j$, $j \in J$, be a system of representatives of the right cosets in $\Gamma_0(p^e n,p)\setminus \Gamma_0(n,p)$ analogous to (\ref{2 system of representatives}).
Then $R_j$, $j \in J$, is already a system of representatives of the right cosets in $\Gamma_0(p^e n)\setminus \Gamma_0(n)$.
\end{itemize}
\end{lemma}

\begin{proof}
\begin{itemize}
\item 
For arbitrary $p,m,n \in \N$ consider a system of representatives $R_j$, $j \in J$ a suitable index set, of the cosets in $\Gamma_0(nm,p) \setminus \Gamma_0(n,p)$.
Let $j_1$, $j_2$ be distinct and consider $R_{j_1} \, R_{j_2}^{-1}$.
Since 
\[
\Gamma_0(nm,p) \, R_{j_1} \cap \Gamma_0(nm,p) \, R_{j_2} = \emptyset 
\]
holds we have
\[
 \Matrix{a}{b}{c}{d} := R_{j_1} \, R_{j_2}^{-1} 
\left\{ \begin{array}{l}
\not\in \Gamma_0(nm,p) \mbox{ and}\\
\in \Gamma_0(n,p).
\end{array} \right.
\]
The entries $b$ and $c$ satisfy $p \mid b$, $n \mid c$ and $nm \not\mid c$.
Hence
\[
R_{j_1} \, R_{j_2}^{-1} 
\left\{ \begin{array}{l}
\not\in \Gamma_0(nm),\\
\in \Gamma_0(n)
\end{array} \right.
\]
which shows that 
\[
\Gamma_0(nm) \, R_{j_1} \cap \Gamma_0(nm) \, R_{j_2} = \emptyset 
\]
\item
For $p$ prime, $p \mid n$ let $\tilde{R}_i$, $i \in I_{p^e n,n}$, be a system of representatives of the right cosets in $\Gamma_0(p^e n)\setminus \Gamma_0(n)$ as defined in (\ref{2 system of representatives}).
We may choose $\tilde{R}_i \in \Gamma_0(n,p)$:
For $\Matrix{a}{b}{c}{d} = \tilde{R}_i$ we know $\gcd(b,d)=1$ and $\gcd(p,d)=1$ since $p \mid c$.
There exists a $t \in \Z$ such that $p \mid b+td$.
Hence $T^t \tilde{R}_i = \Matrix{a+tc}{b+td}{c}{d} \in \Gamma_0(n,p)$ holds.

We find
\begin{eqnarray*}
\bigcup_{i \in I_{p^e n,n}} \Gamma_0(p^e n,p) \, \tilde{R}_i
&=&
\bigcup_{i \in I_{p^e n,n}} \left( \Gamma_0(p^e n) \cap \Gamma_0(n,p) \right)\, \tilde{R}_i\\
&=&
\left( \bigcup_{i \in I_{p^e n,n}}  \Gamma_0(p^e n) \, \tilde{R}_i \right) \cap \Gamma_0(n,p) \\
&=&
\Gamma_0(n) \cap \Gamma_0(n,p)  = \Gamma_0(n,p)
\end{eqnarray*}
and similarly for distinct $i_1$, $i_2 \in I_{p^en,n}$
\begin{eqnarray*}
\emptyset 
&=&
 \left( \Gamma_0(p^e n) \, \tilde{R}_{i_1} \cap  \Gamma_0(p^e n) \, \tilde{R}_{i_1} \right) \cap \Gamma_0(n,p)\\
&=&
 \Gamma_0(p^e n,p) \, \tilde{R}_{i_1} \cap  \Gamma_0(p^e n,p) \, \tilde{R}_{i_1}.
\end{eqnarray*}
Hence $\tilde{R}_i$ is also a system of representatives of $\Gamma_0(p^e n,p)\setminus \Gamma_0(n,p)$.
This shows that the cardinality of the indexsets $I_{p^e n,n}$ and $J$ is equal and $R_j$ is also a system of representatives of the right cosets in $\Gamma_0(p^e n)\setminus \Gamma_0(n)$.
\end{itemize}
\end{proof}

\begin{definition}
\label{B2}
For $g \in \Gnull{nm}$ define $\overline{g}^{(m)} = \overline{g} \in \Gamma_0(n,m)$ by $B_m \,g = \overline{g} \, B_m$.
\end{definition}

Then we can show
\begin{lemma}
\label{B3}
For $p$ prime, $p \mid n$ let $R_j^{p^{e+1}n,pn}$, $j \in I_{p^{e+1}n,pn}$, be a system of representatives of the cosets in $\Gnull{p^{e+1}n}\setminus \Gnull{pn}$ analogous to (\ref{2 system of representatives}).
Then the set $\overline{R_j^{p^{e+1}n,pn}}^{(p)}$, $j \in I_{p^{e+1}n,pn}$, is a system of representatives of the right cosets in $\Gamma_0(p^e n)\setminus \Gamma_0(n)$.
\end{lemma}

\begin{proof}
We have to show that $\overline{R_j^{p^{e+1}n,pn}}=\overline{R_j^{p^{e+1}n,pn}}^{(p)}$, $j \in I_{p^{e+1}n,pn}$, satisfy the property analogous to (\ref{2 system of representatives}).
For this consider $\left( \bigcup_j \Gamma_0(p^e n,p)\, \overline{R_j^{p^{e+1}n,pn}} \right) B_p$.
By Definition~\ref{B2} and Lemma~\ref{B1} we have
\begin{eqnarray*}
\left( \bigcup_j \Gamma_0(p^e n,p)\, \overline{R_j^{p^{e+1}n,pn}} \right) B_p
&=&
\bigcup_j \Gamma_0(p^e n,p)\, B_p \,R_j^{p^{e+1}n,pn}\\
&=&
B_p \, \left(\bigcup_j \Gnull{p^{e+1}n} \,R_j^{p^{e+1}n,pn} \right).
\end{eqnarray*}
Since $R_j^{p^{e+1}n,pn}$, $j \in I_{p^{e+1}n,pn}$, is a system of representatives of the right cosets of $\Gnull{np^{e+1}}$ in $\Gnull{np}$ we find
\[
\left( \bigcup_j \Gamma_0(p^e n,p)\, \overline{R_j^{p^{e+1}n,pn}} \right) B_p = B_p \, \Gnull{pn} = \Gamma_0(n,p) \, B_p.
\]
Similarly we find
\begin{eqnarray*}
&& \!\!\!\!\!\!\!\!\!\!
\Gamma_0(p^e n,p)\, \overline{R_{j_1}^{p^{e+1}n,pn}} \, B_p
\cap
\Gamma_0(p^e n,p)\, \overline{R_{j_2}^{p^{e+1}n,pn}} \, B_p\\
&=&
\Gamma_0(p^e n,p)\, B_p \,R_{j_1}^{p^{e+1}n,pn}
\cap
\Gamma_0(p^e n,p)\, B_p \,R_{j_2}^{p^{e+1}n,pn}\\
&=&
B_p \, \left(\Gnull{p^{e+1}n} \,R_{j_1}^{p^{e+1}n,pn} \cap \Gnull{p^{e+1}n} \,R_{j_2}^{p^{e+1}n,pn} \right)\\
&=&
B_p \, \emptyset = \emptyset \, B_p
\qquad \mbox{for all distinct } j_1,j_2 \in I_{p^{e+1}n,pn}
\end{eqnarray*}
proving the property analogous to (\ref{2 system of representatives}) for $\Gamma_0(p^e n,p)\setminus \Gamma_0(n,p)$.
Lemma~\ref{B8} shows that the set $\overline{R_j^{p^{e+1}n,pn}}^{(p)}$, $j \in I_{p^{e+1}n,pn}$, is a system of representatives of the right cosets in $\Gamma_0(p^e n)\setminus \Gamma_0(n)$.
\end{proof}

Next we show
\begin{lemma}
\label{B5}
For $m,m^\prime$ coprime let $R_j^{mm^\prime n,mn}$, $j \in I_{mm^\prime n,mn}$, be a system of representatives of the cosets in $\Gnull{mm^\prime n}\setminus \Gnull{mn}$.
Then the set $\overline{R_j^{mm^\prime n,mn}}^{(m)}$, $j \in I_{mm^\prime n,m n}$, is a system of representatives of the right cosets in $\Gamma_0(m^\prime n)\setminus\Gamma_0(n)$.
\end{lemma}

\begin{proof}
As in the proof of Lemma~\ref{B3} we have to show that the matrices $\overline{R_j^{mm^\prime n,mn}}=\overline{R_j^{mm^\prime n,mn}}^{(m)}$, $j \in I_{mm^\prime n,mn}$, satisfy the property analogous to (\ref{2 system of representatives}).
For $ \left( \bigcup_j \Gamma_0(m^\prime n,m)\, \overline{R_j^{mm^\prime n,mn}} \right) B_m$
we find according to Definition~\ref{B2} 
\begin{eqnarray*}
\left( \bigcup_j \Gamma_0(m^\prime n,m)\, \overline{R_j^{mm^\prime n,mn}} \right) B_m
&=&
\bigcup_j \Gamma_0(m^\prime n,m)\, B_m \,R_j^{mm^\prime n,mn}\\
&=&
B_m \, \left(\bigcup_j \Gnull{mm^\prime n} \,R_j^{mm^\prime n,mn} \right).
\end{eqnarray*}
Since $R_j^{mm^\prime n,mn}$, $j \in I_{mm^\prime n,m n}$, is a system of representatives of the right cosets in $\Gnull{mm^\prime n}\setminus\Gnull{mn}$ we find
\[
\left( \bigcup_j \Gamma_0(m^\prime n,m)\, \overline{R_j^{mm\prime n,mn}} \right) B_m = B_m \, \Gnull{mn} = \Gamma_0(n,m) \, B_m.
\]
Similarly we find
\begin{eqnarray*}
&& \!\!\!\!\!\!\!\!\!\!
\Gamma_0(m^\prime n,m)\, \overline{R_{j_1}^{mm^\prime n,mn}} \, B_m
\cap
\Gamma_0(m^\prime n,m)\, \overline{R_{j_2}^{mm^\prime n,mn}} \, B_m\\
&=&
\Gamma_0(m^\prime n,m)\, B_m \,R_{j_1}^{mm^\prime n,mn}
\cap
\Gamma_0(m^\prime n,m)\, B_m \,R_{j_2}^{mm^\prime n,mn}\\
&=&
B_m \, \left( \Gnull{mm\prime n} \,R_{j_1}^{mm^\prime n,mn} \cap \Gnull{mm^\prime n} \,R_{j_2}^{mm^\prime n,mn} \right)\\
&=&
B_m \, \emptyset = \emptyset \, B_m= \emptyset
\qquad \mbox{for all distinct } j_1,j_2 \in I_{mm^\prime n,mn}.
\end{eqnarray*}
Hence $\overline{R_{j}^{mm^\prime n,mn}}$, $j \in I_{mm^\prime n,mn}$, is a system of representatives of the right cosets in $\Gamma_0(m^\prime n,m)\setminus \Gamma_0(n,m)$.
The first part of Lemma~\ref{B8} shows that the set of all $\overline{R_{j}^{mm^\prime n,mn}}$ are also a subsystem of representatives of the right cosets in $\Gamma_0(m^\prime n)\setminus \Gamma_0(n)$.
If we show that
\begin{equation}
\label{B5b}
\abs{\Gamma_0(m^\prime n,m)\setminus \Gamma_0(n,m)}
=
\abs{\Gnull{m^\prime n} \setminus \Gnull{n}},
\end{equation}
i.e., $\overline{R_{j^\prime}^{np^{e+1},np}}^{(p)}$, $j^\prime \in I_{np^{e-1},np}$, runs though all cosets in $\Gnull{np^e} \setminus \Gnull{n}$, then set of all $\overline{R_{j}^{mm^\prime n,mn}}$ is already a system of representatives $\Gamma_0(m^\prime n) \setminus \Gamma_0(n)$.

Recall the index formula
\begin{equation}
\label{B5c}
\abs{\Gnull{n} \setminus \Gnull{1}} = n \, \prod_{q \mid n} \left( 1 + \frac{1}{q} \right)
\end{equation}
where the $q$ runs through all prime divisors of $n$.
For the right hand side of (\ref{B5b}) we find
\begin{eqnarray*}
\abs{\Gnull{m^\prime n} \setminus \Gnull{n}}
&=&
\frac{\abs{\Gnull{m^\prime n} \setminus \Gnull{1}}}{\abs{\Gnull{n} \setminus \Gnull{1}}}\\
&=&
\frac{ m^\prime \, \prod_{q \mid m^\prime n} \left( 1 + \frac{1}{q} \right)}{\prod_{r \mid n} \left( 1 + \frac{1}{r} \right)}.
\end{eqnarray*}
Since by construction of $\overline{R_{j}^{mm^\prime n,mn}}$ we have 
\[
\abs{\Gamma_0(m^\prime n,m)\setminus \Gamma_0(n,m)} = \abs{ \Gamma_0(m m^\prime n) \setminus \Gamma_0(mn)}
\]
we find for the left hand side of (\ref{B5b}) 
\begin{eqnarray*}
\abs{\Gamma_0(m^\prime n,m)\setminus \Gamma_0(n,m)}
&=&
\abs{ \Gamma_0(m m^\prime n) \setminus \Gamma_0(mn)}\\
&=&
\frac{\abs{\Gnull{mm^\prime n} \setminus \Gnull{1}}}{\abs{\Gnull{mn} \setminus \Gnull{1}}}\\
&=&
\frac{ mm^\prime n \, \prod_{q \mid mm^\prime n} \left( 1 + \frac{1}{q} \right)}{mn \, \prod_{r \mid mn} \left( 1 + \frac{1}{r} \right)}\\
&=&
\frac{ m^\prime \, \prod_{q \mid m^\prime n} \left( 1 + \frac{1}{q} \right) \, \prod_{q \mid m \atop q \not\mid n} \left( 1 + \frac{1}{q} \right)}{\prod_{r \mid n} \left( 1 + \frac{1}{r} \right) \, \prod_{r \mid m \atop r \not\mid n} \left( 1 + \frac{1}{r} \right)}\\
&=&
\frac{ m^\prime \, \prod_{q \mid m^\prime n} \left( 1 + \frac{1}{q} \right)}{\prod_{r \mid n} \left( 1 + \frac{1}{r} \right)}.
\end{eqnarray*}
Hence Equation~(\ref{B5b}) holds since both sides are equal.
\end{proof}

We need also the following
\begin{lemma}
\label{B4}
For $p$ prime, $p \mid n$ let $R_l^{np^{e+1},n}$, $l \in I_{np^{e+1},n}$, and $R_i^{np,n}$, $i \in I_{np,n}$, denote systems of representatives of the right cosets in $\Gnull{np^{e+1}}\setminus\Gnull{n}$ respectively $\Gnull{np}\setminus\Gnull{n}.$
There exists a system of representatives $R_j^{np^e,n}$, $j \in I_{np^e,n}$ of the right cosets in $\Gnull{np^e}\setminus\Gnull{n}$ such that
\begin{equation}
\label{B4a}
\sum_l B_{p^{e+1}} \, R_l^{np^{e+1},n} = \sum_i \left( \sum_j B_{p^e} \, R_j^{np^e,n} \right) \, B_p \, R_i^{np,n}.
\end{equation}
\end{lemma}

\begin{proof}
Since $[\Gnull{n}:\Gnull{np^{e+1}}]= [\Gnull{n}:\Gnull{np}] \, [\Gnull{np}: \Gnull{np^{e+1}}]$ and therefore $\Gnull{np^{e+1}} \backslash \Gnull{np} \times \Gnull{np} \backslash \Gnull{n} \simeq \Gnull{np^{e+1}} \backslash \Gnull{n}$ we can choose the matrices $R_\cdot^{\cdot,\cdot}$ such that the identity
\[
 R_l^{np^{e+1},n} =   R_{j^\prime}^{np^{e+1},np} \, R_i^{np,n}
\]
holds.
Hence we can identify $l \in I_{np^{e+1},n}$ uniquely with pairs $(j^\prime,i) \in I_{np^{e-1},np} \times I_{np,n}$.
Therefore we find
\[
\sum_l B_{p^{e+1}} \, R_l^{np^{e+1},n} = \sum_{i,j^\prime} B_{p^e}\,B_p \, R_{j^\prime}^{np^{e+1},np} \, R_i^{np,n}.
\]
Definition~\ref{B2} implies 
\[
\sum_l B_{p^{e+1}} \, R_l^{np^{e+1},n} = \sum_{i,j^\prime} B_{p^e} \, \overline{R_{j^\prime}^{np^{e+1},np}}^{(p)} \, B_p \, R_i^{np,n}.
\]
Lemma~\ref{B3} shows that $\overline{R_{j^\prime}^{np^{e+1},np}}^{(p)}$, $j^\prime \in I_{np^{e+1},np}$ is in fact a system of representatives of the right cosets in $\Gnull{np^e}\setminus\Gnull{n}$.
Hence for each $j^\prime$ there exists a unique $j\in I_{np^e,n}$ such that $\Gnull{np^e}\,\overline{R_{j^\prime}^{np^{e+1},np}}^{(p)} = \Gnull{np^e}\, R_j^{np^e,n}$ holds.
This proves Relation (\ref{B4a}).
\end{proof}

Our main result of this section is
\begin{lemma}
\label{B6}
For $m,m^\prime$ coprime let $R_l^{mm^\prime n,n}$, $l \in I_{mm^\prime n,n}$, and $R_i^{mn,n}$, $i \in I_{mn,n}$, denote systems of representatives of the right cosets in $\Gnull{mm^\prime n}\setminus\Gnull{n}$ respectively $\Gnull{mn}\setminus\Gnull{n}$.
There exists a system of representatives $R_j^{m^\prime n,n}$, $j \in I_{m^\prime n,n}$ of the right cosets in $\Gnull{m^\prime n}\setminus\Gnull{n}$ such that
\begin{equation}
\label{B6a}
\sum_l B_{mm^\prime} \, R_l^{mm^\prime n,n} = \sum_i \left( \sum_j B_{m^\prime} \, R_j^{m^\prime n,n} \right) \, B_m \, R_i^{mn,n}.
\end{equation}
\end{lemma}

\begin{proof}
Since $\Gnull{mm^\prime n} \backslash \Gnull{mn} \times \Gnull{mn} \backslash \Gnull{n} \simeq \Gnull{mm^\prime n} \backslash \Gnull{n}$ we can choose the matrices $R_\cdot^{\cdot,\cdot}$ such that the identity
\[
R_l^{mm^\prime n,n} = R_{j^\prime}^{mm^\prime n,mn} \, R_i^{mn,n}
\]
holds.
Hence we can identify $l \in I_{mm^\prime n,n}$ uniquely with pairs $(j^\prime,i) \in I_{m^\prime n,mn} \times I_{mn,n}$.
We find
\[
\sum_l B_{mm^\prime} \, R_l^{mm^\prime n,n} = \sum_{i,j^\prime} B_{m^\prime}\,B_m \, R_{j^\prime}^{mm^\prime n,mn} \, R_i^{mn,n}.
\]
Definition~\ref{B2} implies 
\[
\sum_l B_{mm^\prime} \, R_l^{mm^\prime n,n} = \sum_{i,j^\prime} B_{m^\prime} \, \overline{R_{j^\prime}^{mm^\prime n,mn}} \, B_m \, R_i^{mn,n}.
\]
Lemma~\ref{B5} shows that $\overline{R_{j^\prime}^{mm^\prime n,mn}}$, $j^\prime \in I_{mm^\prime n,mn}$ is indeed a system of representatives of the right cosets in $\Gnull{nm^\prime}\setminus\Gnull{n}$.
Hence for each $j^\prime$ there exists a unique $j\in I_{nm^\prime,n}$ such that $\Gnull{nm^\prime}\,\overline{R_{j^\prime}^{mm^\prime n,mn}} = \Gnull{nm^\prime}\, R_j^{m^\prime n,n}$ holds.
This proves Relation (\ref{B6a}).
\end{proof}

\section{Maass cusp forms and vector valued period functions for $\Gnull{n}$}
\label{periodfunctions}

A \emph{Maass cusp form} $u$ for the congruence subgroup $\Gnull{n}$ is a real-analytic function $u:\HH \to \C$ satisfying:
\begin{enumerate}
\item
$u(gz) =u(z)$ for all $g \in \Gnull{n}$,
\item
$\Delta u= \beta(1-\beta)u$ for some $\beta \in \C$ where $\Delta = -y^2(\partial_x^2+\partial_y^2)$ is the hyperbolic Laplace operator.
We call the parameter $\beta$ the \emph{spectral parameter} of $u$.
\item
$u$ is of rapid decay in all cusps:
if $p \in \Q \cup \{\infty\}$ is a cuspidal point for $\Gnull{n}$ and $g\in\Gmod$ is such that $gp=\infty$ then $u(gz) = \mathrm{O}\left(\im{z}^C\right)$ as $\im{z} \to \infty$ for all $C \in \R$.
\end{enumerate} 
We denote the space of Maass cusp forms for $\Gnull{n}$ with spectral value $\beta$ by $S(n,\beta)$.

\begin{remark}
Item (3) above seems to state two different conditions, namely vanishing in all cusps and an explicit growth condition.
However both these conditions are equivalent since it is shown in \cite{Iw02} that vanishing in the cusp $p=g^{-1} \infty$ implies the stronger growth condition $u(gz) = \mathrm{O}\left(e^{-\im{z}}\right)$ as $\im{z} \to \infty$ which again implies vanishing at the cusp $p=g^{-1}\infty$.
\end{remark}

Consider next the representation $\rho =\rho_\Gnull{n}$ of $\SL{\Z}$ induced from the trivial representation of $\Gnull{n}$ with
\begin{equation}
\label{representation}
\rho(g) = \rho_\Gnull{n}(g) := \big( \delta_\Gnull{n}(R^{n,1}_i\,g\,(R^{n,1}_j)^{-1}) \big)_{1 \leq i,j \leq \mu_n}
\end{equation}
where $R^{n,1}_j$, $j\in I_n=I_{n,1}$, are representatives of the cosets in $\Gnull{n}\setminus \SL{\Z}$ described in (\ref{2 system of representatives}) and
\begin{equation}
\label{delta}
\delta_\Gnull{n}(g) =
\left\{ \begin{array}{ll}
1 \quad &\mbox{if } g \in \Gnull{n} \mbox{ and}\\
0 \quad &\mbox{if } g \not\in \Gnull{n}.
\end{array}\right.
\end{equation}

A \emph{vector valued Maass cusp form} $\vec{u}$ for the congruence subgroup $\Gnull{n}$ with \emph{spectral parameter} $\beta$ is a vector of real-analytic functions $u_i:\HH \to \C$, $i \in I_n$ satisfying 
\begin{enumerate}
\item
$\vec{u}(gz) = \rho(g)\, \vec{u}(z)$ for all $g \in \SL{\Z}$, 
\item
$\Delta u_i= \beta(1-\beta) \, u_i$ for all $i\in I_n$ and
\item
\label{A1.3}
$u_i(z) = \mathrm{O}\left(\im{z}^C\right)$ as $\im{z} \to \infty$ for all $C \in \R$.
\end{enumerate} 
We denote the space of Maass cusp forms for $\Gnull{n}$ with spectral value $\beta$ by $S_\mathrm{ind}(n,\beta)$.

In \cite{Mu05} it is shown that the map $\Pi_nu: S(n,\beta) \to S_\mathrm{ind}(n,\beta)$ given by
\begin{equation}
\label{vvM}
(\Pi_n u)_j(z)=(\vec{u})_j (z)= u\big( R^{n,1}_j z\big)
\qquad  j\in I_n=I_{n,1},
\end{equation}
is bijective.

Denote by $R_\zeta(z)$  the Poisson kernel $R_\zeta(z)= \frac{y}{(x-\zeta)^2 + y^2}$ for $z=x+iy \in \HH$ and $\zeta \in \C \setminus \{x\}$ and by $\eta$ the 1-form
\[
\eta(u,v) = (v \partial_y u -u \partial_y v)dx + (u \partial_x v -v \partial_x u)dy
\]
on $\mathbb{H}$.

It is shown in \cite{Mu05} that the integral transformation $\PP_n: S_\mathrm{ind}(n,\beta) \to \tilde{\FE}(n,\beta)$ with
\begin{equation}
(\PP_n\vec{u})_i(\zeta)=\int_0^{i\infty}\eta(u_i,R_\zeta^\beta)
\end{equation}
maps the space $S_\mathrm{ind}(n,\beta)$ bijectively onto the space $\tilde{\FE}(n,\beta)$ of vector valued  functions $\vec{\Phi}$ for the group $\Gnull{n}$ obeying the equation
\begin{equation}
\label{new three term}
\vec{\phi}(z) - \rho_\Gnull{n}(T^{-1})\, \vec{\phi}(z+1) - (z+1)^{-2\beta} \, \rho_\Gnull{n}\big({T^\prime}^{-1}\big)\, \vec{\phi}\left(\frac{z}{z+1}\right) = 0.
\end{equation}
According to \cite{MM05} we can identify however $\tilde{\FE}(n,\beta)$ with $\FE(n,\beta)$ using the fact that $\chi_{\overline{\Gnull{n}}}$ restricted to $\Gnull{n}$ and $\rho_\Gnull{n}$ are unitarily equivalent.
From now on we therefore replace $\chi_{\overline{\Gnull{n}}}$ when restricted to $\Gnull{n}$ by $\rho_\Gnull{n}$.
Hence vector valued period functions $\vec{\Phi}$ satisfy the three-term Functional Equation (\ref{new three term}).

\smallskip

In the following we use the notation $(f|_0 h)(z)$ for integer matrices with positive determinant for $(f|_0 h)(z)= f(hz)$. The same notation is used also for $\mathbb{Z}$- linear combinations of such matrices.

Since $\Gnull{nm}\subset\Gnull{n}$  any Maass cusp form $u\in S(n,\beta)$ for the group $\Gnull{n}$ is also a Maass cusp form for the group $\Gnull{nm}$ and hence $u\in S(nm,\beta)$.
Denote the corresponding vector valued Maass cusp forms by $\vec{u}_n=\Pi_nu \in S_\mathrm{ind}(n,\beta)$ respectively by $\vec{u}=\Pi_{nm} u\in S_\mathrm{ind}(nm,\beta)$.
It is not difficult to show

\begin{lemma}
\label{n,nm}
The vector valued period functions $\PP_{nm}\vec{u}$  respectively $\PP_{n}\vec{u}_n$  are related by $(\PP_{nm}\vec{u})_j=(\PP_{n}\vec{u}_n)_{\sigma_{m,n}(j)}$ if the representatives $\{R_j^{nm,1}\}_{j\in I_{nm}}$ of the cosets in $\Gnull{nm}\setminus \Gnull{1}$ and the representatives $\{R_i^{n,1}\}_{ i\in I_{n}}$ of the cosets in $\Gnull{n}\setminus \Gnull{1}$ are related by $R_j ^{nm,1}=R_{r(j)}^{nm,n} R_{\sigma_{m,n}(j)}^{n,1},$ where $R_{r(j)} ^{nm,n}$ is a suitable representative of a coset in $\Gnull{nm}\setminus\Gnull{n}$.
\end{lemma}

\begin{proof}
Using (\ref{vvM}) one gets for the representatives $R_\cdot^{mn,1}$ and $R_\cdot^{n,1}$ related as in the Lemma 
\begin{eqnarray*}
\left( \PP_{nm}\vec{u} \right)_j (\zeta)
&=&
\int_0^{i\infty}\eta\left((\Pi_{nm}u)_j ,R_\zeta^\beta \right)
=
\int_0^{i\infty}\eta\left(u\big|_0 R_j^{nm,1} ,R_\zeta^\beta \right)\\
&=&
\int_0^{i\infty}\eta\left(u \big|_0 R_{\sigma_{m,n}(j)}^{n,1} ,R_\zeta^\beta \right)
=
\int_0^{i\infty}\eta\left((\Pi_{n}u)_{\sigma_{m,n}(j)} ,R_\zeta^\beta \right)\\
&=&
\left( \PP_n \vec{u}_n \right)_{\sigma_{m,n}(j)} (\zeta)
\end{eqnarray*}
for all $j \in I_{nm}$.
\end{proof}

\smallskip

For $B_m = \Matrix{m}{0}{0}{1}$ and ${R_j^{nm,n}}$ a system of representatives of the cosets in $\Gnull{nm} \setminus \Gnull{n}$ the following proposition can be easily shown, using similar arguments as in \cite{AL70} for holomorphic modular forms.

\begin{proposition}
\label{2 prop 1}
Let $m$ and $n$ be natural numbers.
\begin{itemize}
\item
If $u$ is in $S(n,\beta)$ then $v=u\big|_0 B_m$ is in $S(nm,\beta)$.
\item
If $u$ is in $S(nm,\beta)$ then $u'=u \big|_0 \sum_{j\in I_{nm,n}} R_j^{nm,n}$ is in $S(n,\beta)$.
\end{itemize}
\end{proposition}

The Maass cusp form $v$ in this Proposition is called an \emph{old Maass cusp form} for $\Gnull{nm}$.
Consider then the vector valued period function $\PP_{nm}\vec{v}$ for $\vec{v}=\Pi_{nm} v$ the vector valued Maass cusp form associated to $v$.
As shown in \cite{MM05} $\PP_{nm}\vec{v}$ can be expressed through the vector valued period function $\PP_{nm}\vec{u}$ with $\vec{u}=\Pi_{nm} u$ for  $u\in S(n,\beta)\subset S(nm,\beta)$ if $v=u|0B_m$ as
\begin{eqnarray}
\label{eq}
&&\big(\PP_{nm}\vec{v}\big)_j(\zeta)=\\
\nonumber
&&\!\!
\sum_{r=1}^{L(\sigma_{R_j^{nm,1}}(B_m)0)}\big[\rho_{\Gnull{nm}}((m_r^{(B_m,j)})^{-1})\PP_{nm}\vec{u}\big]_{\Phi_{B_mṃ}(j)}
\Big|_{\beta}\big(m_r^{(B_m,j)}\sigma_{R_j^{nm,1}}(B_m)\big)(\zeta)
\end{eqnarray}
for all $j\in I_{nm}$.

Thereby the maps $\sigma_{R_j^{nm,1}}:X_{m}\to X_{m}$ and $\Phi_A:I_{nm}\to I_{nm}$ are defined for any representative $R_j^{nm,1}$ of the right cosets in $\Gnull{nm}\setminus \Gnull{1}$ and any $A\in X_{m}$ through the relation
\begin{equation}
\label{sigma Phi}
A \, R_j^{nm,1} \, (\sigma_{R_j^{nm,1}}(A))^{-1} \in \Gnull{nm} \, R_{\Phi_A(j)}^{nm,1}.
\end{equation}
The matrices $m_r$ on the other hand are  defined for any rational number $q\in [0,1)$ as follows (see \cite{MM05} and \cite{Mu05}):
there exists a unique sequence $\{y_0, \ldots ,y_{L(q)}\}$ of rational numbers $y_r= \frac{a_r}{b_r}$ with $\gcd(a_r,b_r)=1$, $b_r \geq 0$ and $y_0=\frac{-1}{0}=-\infty$, $y_1=\frac{0}{1}=0$, $y_{L(q)}=q$, $0<y_r<1$ for $r=2,\ldots, L(q)$ such that $y_0< y_1 < \ldots < y_{L(q)}$ and $\det\Matrix {a_{r-1}}{a_r}{b_{r-1}}{b_r}=-1 $ for $r=1, \ldots, L(q)$.
Then the matrix $m_r$ is given as $m_r=\Matrix{b_r}{-a_r}{b_{r-1}}{-a_{r-1}}$.
Furthermore we can define an element $M(q)$ of the ring $\Z[\SL{\Z}]$ through
\begin{equation}
\label{M(q)}
M(q)=\sum_{r=1}^{L(q)}m_r.
\end{equation}
since obviously $m_r\in \SL{\Z}$ for $1\leq r\leq L(q)$ and $m_1=\Matrix{1}{0}{0}{1}$.
The matrices $\left\{ m_r^{(B_m,j)} \right\}$ in Equation~(\ref{eq}) then correspond to the rational number $q= \sigma_{R_j^{nm,1}}(B_m)(0)$.

It was shown in \cite{MM05} and \cite{Fr05} that $M(q)=\sum_{r=1}^{L(q)} m_r$ and the map $K$ in (\ref{map K 1}) are closely related as follows:
for each $A \in X_n$ and $m_r$ in expression $M(A0)=\sum_{r=1}^{L(A0)} m_r$ we have
\begin{equation}
\label{mr K}
K^l(A) = m_{l+1}A
\qquad \mbox{for all } l \in \{1,\ldots,L(A0) -1\}
\end{equation}
and in particular  $k_A =L(A0)-1$.

Hence Lemma~\ref{lemma 1} is equivalent to
\begin{equation}
\label{lemma 2a}
A_{l_{i,j}} \, \big( K^j (A_{\sigma_{n,m}(i)}) \big) \,A_i^{-1}
=
A_{l_{i,j}} \, m_{j+1} \, A_{\sigma_{n,m}(i)} \, A_i^{-1}\in SL(2,\mathbb{Z}),
\end{equation}
where the $m_j$'s are given by (\ref{M(q)}) with $q=A_{\sigma_{n,m}(i)}0$.

\smallskip

Consider next a Maass cusp form $u\in S(nm,\beta)$ respectively its vector valued form $\vec{u}= (u_j)_{j\in I_{nm}} =\Pi_{nm} u$ with entries $u_j=u\big|_0 R_j^{nm,1}$.
Then one gets

\begin{lemma}
\label{Lem}
The vector valued period function $\PP_n(\Pi_n v)$ of the Maass cusp form $v = \sum_{i\in I_{mn,n}} u\big|_0 R_i^{nm,n} \in S(n,\beta)$ is related to the vector valued period function $\PP_{nm}\vec{u}$ of the Maass cusp form $u\in S(nm,\beta)$ by $(\PP_n\vec{v})_j = \sum_{l\in\sigma_{m,n}^{-1}(j)} (\PP_{nm}\vec{u})_l$, $j \in I_{nm}$.
\end{lemma}

\begin{proof}
By Proposition~\ref{2 prop 1} the function $v = \sum_{i \in I_{mn,n}} u\big |_0 R_i^{nm,n}$, with $I_{mn,n}$ the index set of the right cosets of $\Gnull{nm}$ in $\Gnull{n}$, is indeed in $S(n,\beta)$.
Hence $\Pi_n v =\left(v\big|_0 R_j^{n,1}\right)_{j \in I_n}$, where $\{R_j^{n,1}\}_{j \in I_n}$ is a system of representatives of the cosets in $\Gnull{n}\backslash \SL{\Z}$.
But the family $\{ R_i^{nm,n} R_j^{n,1}\}_{i\in I_{mn,n}, j \in I_n}$ is a system of representatives of the cosets in  $\Gnull{nm} \backslash \SL{\Z}$ and hence
\begin{eqnarray*}
(\PP_n \Pi_n v)_j
&=&
\int_0^{i\infty}\eta(v \big|_0 R_j^{n,1},R_\zeta^\beta)
=
\sum_{i \in I_{mn,m}} \int_0^{i \infty} \eta \big( u\big|_0 R_i^{nm,n}R_j^{n,1},R_\zeta^\beta \big)\\
&=&
\sum_{l \in I_{nm} \atop \sigma_{m,n}(l) =j} \int_0^{i \infty} \eta \big( u\big|_0 R_l^{mn,1},R_\zeta^\beta \big)
=
\sum_{l \in\sigma_{m,n}^{-1}(j)} \big( \PP_{nm}\vec{u} \big)_l.
\end{eqnarray*}
\end{proof}

Lemma~\ref{Lem} gives a simple interpretation of the solution $\vec{\phi}$ of Lewis' equation for $\Gnull{n}$ in the first part of Proposition~\ref{prop 1}: it corresponds to a Maass cusp form for $\Gnull{n}$ constructed from such a form for the subgroup $\Gnull{nm}$ as described in the second part of Proposition~\ref{2 prop 1}.

\smallskip

In the following we need
\begin{lemma}
\label{uniqueness}
If for $i \in I_{nm}$  there exist $j,j^\prime \in I_n$, $k,k^\prime \in I_m$ and matrices $\gamma, \gamma^\prime \in \SL{\Z}$ such that $\gamma \,A_j\,A_k = A_i = \gamma^\prime \,A_{j^\prime}\,A_{k^\prime}$,
then $j=j^\prime$, $k=k^\prime$ and $\gamma=\gamma^\prime$.
\end{lemma}

\begin{lemma}
\label{uniqueness 1}
Let $(A_1,B_1)$, $(A_2,B_2) \in X_n^\star \times X_m^\star$ be two different pairs of matrices such that $A_1 B_1= T^k A_2B_2$ holds for some $k \in\Z$.
Then, the matrix entries of $A_1 B_1$ have a common divisor strictly larger than $1$.
\end{lemma}

\begin{proof}
For $A_l=\Matrix{c_l}{b_l}{0}{\frac{n}{c_l}}$, $B_l=\Matrix{d_l}{f_l}{0}{\frac{m}{d_l}}$ we find
\begin{equation}
\label{CAB}
C:=\Matrix{c_1d_1}{c_1f_1+b_1\frac{m}{d_1}}{0}{\frac{nm}{c_1d_1}}=\Matrix{c_2d_2}{c_2f_2+b_2 \frac{m}{d_2}+k\frac{nm}{c_2d_2}}{0}{\frac{nm}{c_2d_2}}.
\end{equation}
Therefore 
\begin{equation}
\label{c_id_i}
c_1d_1=c_2d_2 \quad\mbox{and}\quad c_1f_1+b_1\frac{m}{d_1}=c_2f_2+b_2 \frac{m}{d_2}+k\frac{nm}{c_2d_2}.
\end{equation}
We consider the following cases:
\begin{itemize}
\item
If $d_1=d_2$ then $c_1=c_2$ and hence $c_1(f_1-f_2)=\frac{m}{d_1}(b_2-b_1+k\frac{n}{c_1})$.
Since $|f_1-f_2| \leq\frac{m}{d_1}-1$ it follows that $(c_1,\frac{m}{d_1})>1$.
Indeed, if $c_1=1$ then $b_2-b_1-k\frac{n}{c_1}=0$.
But $|b_1-b_2| \leq\frac{n}{c_1}-1$ implies $k=0$ and hence $b_1=b_2$.
Then also $f_1=f_2$ and therefore $A_1=A_2$ and $B_1=B_2$.
If on the other hand $\frac{m}{d_1}=1$ then $f_1=f_2=0$ and $b_2-b_1-k\frac{n}{c_1}=0$ must hold.
But this is only possible if $b_1=b_2$ and hence we find again $A_1=A_2$ and $B_1=B_2$.
If therefore $(c_1,\frac{m}{d_1})>1$ we get $\gcd(c_1d_1,c_1f_1+b_1\frac{m}{d_1},\frac{n}{c_1}\frac{m}{d_1})>1$.
\item
Consider next the case $(d_1,d_2)=1$.
From $c_1d_1=c_2d_2$ we conclude that $d_2\mid c_1$ and $d_1\mid c_2$.
From this it follows that $d_1\mid\frac{m}{d_2}$.
Otherwise obviously $d_1\mid d_2$.
If $d_1=1$ we can assume $d_2>1$ since the case $d_1=d_2$ has been treated already before.
Then $c_1=c_2d_2$ and therefore $d_2\mid c_1$.
Since also $d_2\mid m$ we get $\gcd(c_1d_1,c_1f_1+b_1m,\frac{n}{c_1}m) >1$.
\item
It remains to consider the case $(d_1,d_2)>1$.
Then there exists a number $q_1>1$ with $d_1=q_1d_1^{(1)}$, $d_2=q_1d_2^{(1)}$ and $m=q_1 m^{(1)}$.
The matrix $C$ therefore has the form
\begin{equation}
C
=
\left( \begin{array}{cc}
c_1 q_1 d_1^{(1)} & c_1 f_1 + b_1 \frac{m^{(1)}}{d_1^{(1)}}\\
0                 & \frac{nm^{(1)}}{c_1 d_1^{(1)}}
\end{array} \right)
=
\left( \begin{array}{cc}
c_2 q_1 d_2^{(1)} & c_2 f_2 + b_2\frac{m^{(1)}}{d_2^{(1)}} + k \frac{nm^{(1)}}{c_2d_2^{(1)}}\\
0                 & \frac{nm^{(1)}}{c_2d_2^{(1)}}
\end{array} \right).
\end{equation}
But these two representations of the matrix $C$ are similar to the ones in equation \ref{CAB}.
Equations  \ref{c_id_i} now read
\begin{equation}
c_1d_1^{(1)}=c_2d_2^{(1)} \quad\mbox{and}\quad c_1f_1+b_1\frac{m^{(1)}}{d_1^{(1)}}=c_2f_2+b_2 \frac{m^{(1)}}{d_2^{(1)}}+k\frac{nm^{(1)}}{c_2d_2^{(1)}}.
\end{equation}
with $m^{(1)}< m$ and $d_l^{(1)}<d_l$ for $l=1,2$.
We can therefore apply our chain of arguments to $m^{(1)}$ and $d_l^{(1)}$, $l=1,2$ to arrive at a new triple $m^{(2)}$ and $d_l^{(2)}$, $l=1,2$ and so on.
Since $m$ is finite we arrive after finitely many steps at the case $(d_1^{r},d_2^{(r)})=1$ or $d_1^{r}=d_2^{(r)}$ which was handled already before.
\end{itemize}
This concludes the proof of Lemma \ref{uniqueness 1}.
\end{proof}

\begin{proof}[Proof of Lemma~\ref{uniqueness}]
Since all $A$'s are upper triangula matrices we see that both $\gamma$ and $\gamma^\prime$ are translation matrices.
We can find a $k \in \Z$ such that $\gamma^{-1} \gamma^\prime = T^k$.

We assume that $(j,k) \not= (j^\prime,k^\prime)$.
Since $C := A_j\, A_k = T^k \, A_{j^\prime} \, A_{k^\prime}$ holds Lemma~\ref{uniqueness 1} implies that the matrix entries of $C$ have common divisor strictly larger than $1$.
Hence the matrix entries of $\gamma \,C$ have common divisor strictly larger than $1$.
But this contradicts the fact $\gamma \, C =A_i \in X_{nm}^\star$.
Hence $(j,k)=(j^\prime, k^\prime)$, proving the lemma.
\end{proof}

\smallskip

Next we want to show that also the solutions $\vec{\Phi}=(\Phi_j)_{j\in I_{nm}}$ in Equation~(\ref{functions Phi}) have a simple interpretation similar to Lemma \ref{Lem}.
But before doing this we have to recall some more notations and facts from \cite{MM05}.

Let us start with the map $h_n:I_n \to I_n$ introduced in \cite{MM05}:
For $i \in I_n$ consider the representative $R_i^{n,1} = \Matrix{a}{b}{c}{d}$ and the matrix $\Matrix{c}{b}{0}{\frac{n}{c}} \in X^\star_n$ with $b \equiv d \bmod \frac{n}{c}$.
By (\ref{matrix Ai}) there exists an index $l(i) \in I_n$ such that $A_{l(i)} = \Matrix{c}{b}{0}{\frac{n}{c}}$.
Then it was shown in \cite{MM05} that the relation
\begin{equation}
\label{Ri Aiota}
\Matrix{0}{1}{-n}{0}\, R_i^{n,1} \in \SL{\Z} \, A_{l(i)}
\end{equation}
holds.
We define a map $h_n:I_n\to I_n$ by
\begin{equation}
\label{def of hn}
h_n(i):= l(i).
\end{equation}
It is shown in \cite{MM05} that $h_n:I_n\to I_n$ is then bijective.

Lateron we need also the following result \cite{Fr05}:
\begin{lemma}
\label{lemma Fr}
For $n \in \N$ and $i \in I_n$ we have $\sigma_{R_i^{n,1}}(B_n) = A_{h_n(i)}$.
\end{lemma}

\begin{proof}
Multiplying (\ref{Ri Aiota}) by $S^{-1}$ from the left, we find
\[
B_n \, R_i^{n,1} \in \SL{\Z} \, A_{h_n(i)}.
\]
On the other hand using (\ref{sigma Phi}) we find
\[
B_n \, R_i^{n,1} \in \SL{\Z} \, \sigma_{R_i^{n,1}}(B_n) 
\]
The fact that both matrices $\sigma_{R_i^{n,1}}(B_n)$ and $A_{h_n(i)}$ are in $X_n^\star$ implies that
\[
\sigma_{R_i^{n,1}}(B_n) = A_{h_n(i)}.
\]
\end{proof}

\begin{lemma}
\label{h sigma}
For $n,m \in N$ and $i \in I_{nm}$ we have $h_n\big(\sigma_{m,n}(i)\big) = \sigma_{n,m} \big( h_{nm}(i) \big)$.
\end{lemma}

\begin{proof}
Let $\{R_i^{nm,1}, i \in I_{nm}\}$, $\{R_l^{mn,n},l \in I_{mn,n}\}$, and $\{R_\nu^{n,1}, \nu \in I_n\}$, be systems of representatives of the right cosets in $\Gnull{mn} \backslash \SL{\Z}$, $\Gnull{mn} \backslash \Gnull{n}$ and $\Gnull{n} \backslash \SL{\Z}$ respectively, such that $R_i^{nm,1} = R_l^{nm,n} \, R_\nu^{n,1}$ holds for all $i \in I_{mn}$ with $\nu := \sigma_{m,n}(i)$ and suitabe $l=l_i \in I_{mn,n}$.

For  $i \in I_{mn}$ consider the matrix $A_{h_{nm}(i)}$.
According to (\ref{Ri Aiota}) there exists a $g \in \SL{\Z}$ such that $g \, A_{h_{nm}}(i) = \Matrix{0}{-1}{nm}{0} \, R_i^{nm,1}$.
Hence using Definition~\ref{B2} we find
\begin{eqnarray*}
g \, A_{h_{nm}}(i)
&=&
\Matrix{0}{-1}{nm}{0} \, R_i^{nm,1}
=
\Matrix{0}{-1}{m}{0} \, \Matrix{n}{0}{0}{1} \, R_l^{nm,n} \, R_\nu^{n,1}\\
&=&
\Matrix{0}{-1}{m}{0} \, \overline{R_l^{nm,n}}^{(1)} \, \Matrix{n}{0}{0}{1} \, R_\nu^{n,1}\\
&=&
\Matrix{0}{-1}{m}{0} \, \overline{R_l^{nm,n}}^{(1)} \, S^{-1} \, \Matrix{0}{-1}{n}{0} \, R_\nu^{n,1}.
\end{eqnarray*}
Again by (\ref{Ri Aiota}) there exists  $g_\nu \in \SL{\Z}$ such that
\[
g \, A_{h_{nm}}(i)
=
\Matrix{0}{-1}{m}{0} \, \overline{R_l^{nm,n}}^{(n)} \,S^{-1} \, g_\nu \, A_{h_n(\nu)}.
\]
Since $\overline{R_l^{nm,n}}^{(1)} \,S^{-1} \, g_\nu$ is a representative of a suitable right coset $\Gnull{m} \, R_{\bar{l}}^{m,1}$ for some $\bar{l} \in I_m$ we find
\[
g \, A_{h_{nm}}(i)
=
\Matrix{0}{-1}{m}{0} \,\gamma_0 \, R_{\bar{l}}^{m,1} \, A_{h_n(\nu)}
=\tilde{\gamma}_0 \, \Matrix{0}{-1}{m}{0} \, R_{\bar{l}}^{m,1} \, A_{h_n(\nu)} \]
 for some $\gamma_0, \tilde{\gamma}_0 \in \Gamma_0(m).$

By (\ref{Ri Aiota}) there exists a $g_{\bar{l}} \in \SL{\Z}$ such that
\begin{equation}
\label{h sigma 1}
g \, A_{h_{nm}}(i)
=
g_{\bar{l}} \, A_{h_m(\bar{l})} \, A_{h_n(\nu)}.
\end{equation}
On the other hand  Lemma~\ref{lemma 1} shows that
\begin{equation}
\label{h sigma 2}
A_{l_{h_{nm}(i),0}} \, A_{\sigma_{m,n}(h_{nm}(i))} \in \SL{\Z}\, A_{h_{nm}(i)}
\end{equation}
Lemma~\ref{uniqueness} then implies that the two factorizations of $A_{h_{nm}(i)}$ in (\ref{h sigma 1}) and (\ref{h sigma 2}) are identical and hence
 $A_{\sigma_{m,n}(h_{nm}(i))} = A_{h_n(\nu)}$ implying
\[
\sigma_{m,n}\big(h_{nm}(i)\big) = h_n \big( \sigma_{m,n}(i) \big)
\]
since $\nu = \sigma_{m,n}((i))$.
\end{proof}

We next recall the representation $\tilde{\rho}_n$ of $\SL{\Z}$ introduced in \cite{MM05} with
\[
\big(\tilde{\rho}_n(g)\big)_{i,j}=\delta_{\SL{\Z}}(A_igA_j^{-1}) \qquad \mbox{for } i,j\in I_n
\]
and $A_i$ as defined in (\ref{matrix Ai}).
In \cite{MM05} it is shown that for any $g \in \SL{\Z}$ one has
\begin{equation}
\label{hn}
R_i^{n,1}g (R_j^{n,1})^{-1}\in \Gnull{n}\quad \mbox{iff} \quad A_{h_n(i)} g A_{h_n(j)}^{-1}\in \SL{\Z}.
\end{equation}
From this one concludes \cite{MM05}
\begin{lemma}
\label{representation roh and tilde(rho)}
The two representations $\rho_{\Gnull{n}}$ and $\tilde{\rho}_n$ are unitarily equivalent with
\begin{equation}
\label{Hn}
\rho_{\Gnull{n}}(g)= \mathcal{H}_n^{-1}\tilde{\rho}_n(g)\mathcal{H}_n
\end{equation}
where $\mathcal{H}_n$ is the $\mu_n \times \mu_n$ matrix with entries
\begin{equation}
(\mathcal{H}_n)_{i,j}=1 \quad \mbox{if} \quad h_n(j)=i \quad \mbox{and} \qquad (\mathcal{H}_n)_{i,j}=0 \quad \mbox{otherwise}.
\end{equation}
\end{lemma}
In terms of the matrix elements Relation~(\ref{Hn}) reads
\begin{equation}
\label{hn 2}
(\rho_{\Gnull{n}}(g))_{i,j}=(\tilde{\rho}_n(g))_{h_n(i),h_n(j)}, \quad i,j \in I_n.
\end{equation}

There is a simple relation between the representations  $\rho_{\Gnull{n}}$ and $\rho_{\Gnull{nm}}$:
\begin{lemma}
\label{rhorhotilde}
For $g\in \SL{\Z}$ the matrix elements $\big(\rho_{\Gnull{n}}(g)\big)_{k,l}$, $k,l\in I_n$, and $\big(\rho_{\Gnull{nm}}(g)\big)_{i,j}$, $i,j\in I_{nm}$ are related by
\[
\big(\rho_{\Gnull{n}}(g)\big)_{k,l}= \sum_{j\in\sigma_{m,n}^{-1}(l)}\big(\rho_{\Gnull{nm}}(g)\big)_{i,j}
\]
for any $i\in\sigma_{m,n}^{-1}(k)$.
\end{lemma}

\begin{proof}
For $(\rho_{\Gnull{nm}}(g))_{i,j}=1$ we have $R_i^{nm,1}g (R_j^{nm,1})^{-1}\in\Gnull{nm}\subset \Gnull{n}$.
The representatives $\{R_i^{nm,1}\}$ of the cosets in $\Gnull{nm}\setminus \SL{\Z}$ however can be written as $R_i^{nm,1}= R_{r(i)}^{nm,n}R_{\sigma_{m,n}(i)}^{n,1}$ with $R_{r(i)}^{nm,n}$ some representative of a coset in  $\Gnull{nm}\setminus \Gnull{n}$ and $ R_{\sigma_{m,n}(i)}^{n,1}$ a representative in $\Gnull{n}\setminus \SL{\Z}$.
Since $R_{r(i)}^{nm,n}\in \Gnull{n}$ we conclude $R_{\sigma_{m,n}(i)}^{n,1}g(R_{\sigma_{m,n}(j)}^{n,1})^{-1}\in \Gnull{n}$ and hence
\[
\big(\rho_{\Gnull{n}}(g)\big)_{\sigma_{m,n}(i),\sigma_{m,n}(j)}=1.
\]

On the other hand assume $(\rho_{\Gnull{n}}(g))_{k,l}=1$ and hence $R_k^{n,1}g(R_l^{n,1})^{-1}\in \Gnull{n}$.
Take any $i \in \sigma_{m,n}^{-1}(k)$.
Then there exists a unique $j\in \sigma_{m,n}^{-1}(l)$ such that $R_i^{nm,1}g (R_j^{nm,1})^{-1} \in\Gnull{nm}$.
Indeed, for the representative $R_i^{nm,1}$ with  $R_i^{nm,1}=R_{r(i)}^{nm,n}R_k^{n,1}$ we find
$R_{r(i)}^{nm,n}R_k^{n,1}g(R_l^{n,1})^{-1}\in \Gnull{n}$.
Therefore there exists a unique $\tilde{r} \in I_{m,n}$ such that $R_{r(i)}^{nm,n}R_k^{n,1}g(R_l^{n,1})^{-1} (R_{\tilde{r}}^{nm,n})^{-1} \in \Gnull{nm}$.
For $R_j^{nm,1}=R_{\tilde{r}}^{nm,n}R_l^{n,1}$ with $\sigma_{m,n}(j)=l$ we then get $R_i^{nm,1} g (R_j^{nm,1})^{-1} \in \Gnull{nm}$.
Obviously this $j$ is uniquely defined for every $i\in I_n$ with $\sigma_{m,n}(i)=k$.
But this concludes the proof of the lemma.
\end{proof}

To compare the vector valued period function $\PP_{nm}\vec{v}$ in (\ref{eq}) and the solution $\vec{\Phi}$ of the Lewis equation for the group $\Gnull{nm}$ in (\ref{functions Phi}) we need some further results from \cite{MM05}:
\begin{lemma}
\label{rhotilde}
For $j\in I_{nm}$ and $M(A_{\sigma_{n,m}(j)}0)=\sum_{s=1}^{L(A_{\sigma_{n,m}(j)}0)} m_s^{(j)}$ one finds
\[
\Big[\tilde{\rho}_n\big((m_s^{(j)})^{-1}\big)\Big]_{l_{j,0},l_{j,r}}=\delta_{s,r+1}
\]
where $\delta_{s,r+1}$ denotes the Kronecker Delta function.
\end{lemma}
Since $L(A_{\sigma_{n,m}(j)}0) = k_{\sigma_{n,m}(j)} +1$ we can write $\phi_{l_{j,s}}$ in expression (\ref{functions Phi}) as
\begin{eqnarray*}
\phi_{l_{j,s}}
&=&
\sum_{r=0}^{k_{\sigma_{n,m}(j)}} \delta_{s+1,r+1} \, \phi_{l_{j,r}}
=
\sum_{r=0}^{k_{\sigma_{n,m}(j)}}[\tilde{\rho}_n((m_{s+1}^{(j)})^{-1})]_{l_{j,0},l_{j,r}} \, \phi_{l_{j,r}}\\
&=&
\Big[\tilde{\rho}_n\big((m_{s+1}^{(j)})^{-1}\big) \, \vec{\phi}\Big]_{l_{j,0}}.
\end{eqnarray*}
Replacing the index $j\in I_{nm}$ by the index $h_{nm}(i)$ with $h_{nm}:I_{nm}\to I_{nm}$ defined in analogy to the map $h_n$ in (\ref{hn}) we find for $\Phi_{h_{nm}(i)}$ in (\ref{functions Phi}) 
\begin{equation}
\label{phihnm1}
\Phi_{h_{nm}(i)}(\zeta)
=
\sum_{j=0}^{k_{\sigma_{n,m}(h_{nm}(i))}} \Big[\tilde{\rho}_n \big( (m_{j+1}^{(h_{nm} (i))})^{-1} \big) \Big] \, \vec{\phi}]_{l_{h_{nm}(i),0}}\big|_\beta K^j \big(A_{\sigma_{n,m}(h_{nm}(i))}\big) (\zeta),
\end{equation}
with $k_{\sigma_{n,m}(i)}=k_{A_{\sigma_{n,m}(i)}}$.
Using (\ref{mr K}) the component $\Phi_{h_{nm}(i)}$ of $\vec{\Phi}$ can be written as
\[
\Phi_{h_{nm}(i)}
=
\sum_{j=0}^{k_{\sigma_{n,m}(i)}} [\tilde{\rho}_n((m_{j+1}^{(h_{nm}(i))})^{-1})\vec{\phi}]_{l_{h_{nm}(i),0}}\big|_\beta m_{j+1}^{(h_{nm}(i))} A_{\sigma_{n,m}\big(h_{nm}(i)\big)}.
\]
Since $k_{\sigma_{n,m}(i)} = L(A_{\sigma_{n,m}(i)}0) -1$ we write this equation as
\begin{equation}
\label{phihnm}
\Phi_{h_{nm}(i)}
=
\sum_{j=1}^{L(A_{\sigma_{n,m}(i)}0)} [\tilde{\rho}_n((m_j^{(h_{nm}(i))})^{-1})\vec{\phi}]_{l_{h_{nm}(i),0}} \big|_\beta m_j^{(h_{nm}(i))} A_{\sigma_{n,m}\big(h_{nm}(i)\big)}.
\end{equation}

On the other hand consider the vector valued period function $\PP_{nm}\vec{v}$ in (\ref{eq}).
It follows from Lemma \ref{rhorhotilde} that $\big(\rho_{\Gnull{nm}}(g)\big)_{i,j}=1$ leads to
\[
\big(\rho_{\Gnull{n}}(g)\big)_{\sigma_{m,n}(i),\sigma_{m,n}(j)}=1
\]
and hence with Lemma \ref{n,nm} one concludes that
\[
\big[\rho_{\Gnull{nm}}(g)\,\PP_{nm}\vec{u}\big]_j
=
\big[\rho_{\Gnull{n}}(g) \, \PP_{n}\vec{u}_n\big]_{\sigma_{m,n}(j)}.
\]
Hence $\PP_{nm}\vec{v}$ in (\ref{eq}) has the form
\begin{eqnarray}
\label{Pnm}
&&
\big(\PP_{nm}\vec{v}\big)_i =\\
\nonumber
&&
\sum_{r=1}^{L(\sigma_{R_j^{nm,1}}(B_m)0)}[\rho_{\Gnull{n}}((m_r^{(B_m,i)})^{-1})\PP_{n}\vec{u}_n]_{\sigma_{m,n}(\Phi_{B_m}(i))} \big|_{\beta} m_r^{(B_m,i)} \sigma_{R_i^{nm,1}}(B_m).
\end{eqnarray}

To relate now $\Phi_{h_{nm}(i)}$ in (\ref{phihnm}) and $(\PP_{nm}\vec{v})_i$ in (\ref{Pnm}) we use \cite{Fr05}
\begin{lemma}
\label{lA}
For $B_m=\Matrix{m}{0}{0}{1}$ the set
\[
\Big\{\big(l_{h_{nm}(i),0}, A_{\sigma_{n,m}(h_{nm}(i))} \big) \Big\}_{i \in I_{nm}}
\]
coincides with the set
\[
\Big\{\Big( h_n\big(\sigma_{m,n}(\Phi_{B_m}(i))\big),\sigma_{R_i^{nm,1}}(B_m)\Big)\Big\}_{i \in I_{nm}}.
\]
Indeed  $l_{h_{nm}(i),0}=h_n\Big(\sigma_{m,n}\big(\Phi_{B_m}(i)\big)\Big)$ and $A_{\sigma_{n,m}(h_{nm}(i))}=\sigma_{R_i^{nm,1}}(B_m)$.
\end{lemma}

To prove Lemma~\ref{lA} we need the following result \cite{Fr05}
\begin{lemma}
\label{lemma2}
For $m,n \in \N$ and  $\{R_i^{mn,1}\}_{i \in I_{mn}}$ respectively $\{R_k^{n,1}\}_{k \in I_n}$ systems of representatives for the cosets in $\Gnull{nm}\backslash \SL{\Z} $ respectively $\{R_k^{n,1}\}_{k \in I_n}$ in $\Gnull{n}\backslash \SL{\Z}$ with $R_i^{mn,1} \in \Gnull{n} \, R_{\sigma_{m,n}}(i)^{n,1}$ and
$B_n=\Matrix{n}{0}{0}{1}$ and $B_m=\Matrix{m}{0}{0}{1}$
one has for all $i \in I_{nm}$
\[
B_n \, B_m \, R_i^{nm,1} \in \SL{\Z} \, A_{h_n(\sigma_{m,n}(\Phi_{B_m}(i)))} \, \sigma_{R_i^{nm,1}}(B_m).
\]
\end{lemma}

\begin{proof}
Take an $i \in I_{nm}$.
From (\ref{sigma Phi}) it follows that 
\[
B_m \, R_i^{nm,1} \in \Gnull{nm} \, R_{\Phi_{B_m}(i)}^{nm,1} \, \sigma_{R_i^{nm,1}}(B_m).
\]
But by assumption
\[
R_{\Phi_{B_m}(i)}^{nm,1} \in \Gnull{n} \, R_{\sigma_{m,n}(\Phi_{B_m}(i))}^{n,1}
\]
and again by (\ref{sigma Phi}) we get
\[
B_n \, R_{\sigma_{m,n}(\Phi_{B_m}(i))}^{n,1} \in \SL{\Z}\, \sigma_{R_{\sigma_{m,n}(\Phi_{B_m}(i))}^{n,1}}(B_n).
\]
Therefore
\[
R_{\Phi_{B_m}(i)}^{nm,1}
\in
\Gnull{n} \,B_n^{-1} \,\SL{\Z} \, \sigma_{R_{\sigma_{m,n}(\Phi_{B_m}(i))}^{n,1}}(B_n)
\]
and hence
\[
B_m \, R_i^{nm,1}
\in
\Gnull{n} \, B_n^{-1} \, \SL{\Z} \, \sigma_{R_{\sigma_{m,n}(\Phi_{B_m}(i))}^{n,1}}(B_n) \, \sigma_{R_i^{nm,1}}(B_m).
\]
Therefore there exists $\gamma_0\in \Gnull{n}$ such that
\[
B_m \, R_i^{nm,1} \in \gamma_0 \, B_n^{-1} \,  \SL{\Z} \, \sigma_{R_{\sigma_{m,n}(\Phi_{B_m}(i))}^{n,1}}(B_n) \, \sigma_{R_i^{nm,1}}(B_m).
\]
A simple calculation shows that $\gamma_0 \,B_n^{-1}= B_n^{-1} \, \gamma^0$ with $\gamma^0 \in \Gamma_0(1,n)\subset \SL{\Z}$.
Hence
\[
B_m \, R_i^{nm,1} \in  B_n^{-1} \, \SL{\Z} \, \sigma_{R_{\sigma_{m,n}(\Phi_{B_m}(i))}^{n,1}}(B_n) \, \sigma_{R_i^{nm,1}}(B_m).
\]
Since Lemma~\ref{lemma Fr} implies
\[
\sigma_{R_{\sigma_{m,n}(\Phi_{B_m}(i))}^{n,1}}(B_n)=A_{h_n(\sigma_{m,n}(\Phi_{B_m}(i)))},
\]
one arrives finally at
\[
B_n \, B_m \, R_i^{nm,1}
\in
\SL{\Z} \, A_{h_n(\sigma_{m,n}(\Phi_{B_m}(i)))} \, \sigma_{R_i^{nm,1}}(B_m).
\]
\end{proof}

\begin{proof}[Proof of Lemma~\ref{lA}]
According to Relation \ref{Ri Aiota}
\[
\Matrix{0}{1}{-nm}{0} \, R_i^{nm,1} \in \SL{\Z} \, A_{h_{nm}(i)}
\]
and therefore also 
\begin{equation}
\label{BnBm}
B_n \, B_m \, R_i^{nm,1} \in \SL{\Z} \, A_{h_{nm}(i)}.
\end{equation} 
Lemma \ref{lemma 1} for $j=0$ shows the existence of $g \in \SL{\Z}$ such that
\[
g \, A_{l_{h_{nm}(i),0}} \, A_{\sigma_{n,m}(h_{nm}(i))} = A_{h_{nm}(i)}.
\]
Inserting this into (\ref{BnBm}) leads to
\[
B_n \, B_m \, R_i^{nm,1} \in \SL{\Z} \, A_{l_{h_{nm}(i),0}} \, A_{\sigma_{n,m}(h_{nm}(i))}.
\]

On the other hand Lemma~\ref{lemma2} shows that
\[
B_n \, B_m \, R_i^{nm,1} \in \SL{\Z} \, A_{h_n(\sigma_{m,n}(\Phi_{B_m}(i)))} \, \sigma_{R_i^{nm,1}}(B_m).
\]
Hence there exists $\gamma \in \SL{\Z}$ such that
\[
g \,A_{h_n(\sigma_{m,n}(\Phi_{B_m}(i)))} \, \sigma_{R_i^{nm,1}}(B_m) = A_{h_{nm}(i)} = \gamma \, A_{l_{h_{nm}(i),0}} \, A_{\sigma_{n,m}(h_{nm}(i))}.
\]
Then Lemma~\ref{uniqueness} implies that
\[
l_{h_{nm}(i),0} = h_n(\sigma_{m,n}(\Phi_{B_m}(i)))
\qquad \mbox{and} \qquad
A_{\sigma_{n,m}(h_{nm}(i))} = \sigma_{R_i^{nm,1}}(B_m).
\]
\end{proof}

Summarizing we therefore have shown \cite{Fr05}
\begin{lemma}
\label{Phi}
Let $\vec{\phi}$ be a solution of the Lewis Equation (\ref{three-term 3}) for $\Gnull{n}$ with $\vec{\phi}=\mathcal{H}_n P_n\vec{u}$ and $\vec{u}=\Pi_n u$ the vector valued Maass form of $u\in S(n,\beta)$. 
Then the solution $\vec{\Phi}=(\Phi_j)_{j\in I_{nm}}$  of (\ref{three-term 3}) for $\Gnull{nm}$  in Theorem \ref{theorem 2} can be expressed as 
\[
\vec{\Phi}= \mathcal{H}_{nm}\,\PP_{nm} \vec{v}.
\]
Thereby 
$\PP_{nm}\vec{v}$ is the old period function in (\ref{eq}) with $\vec{v}= \Pi_{nm}v$ and $v$ the old Maass cusp form $v=u \big|_0 B_m \in S(nm,\beta)$ determined by $u \in S(n,\beta)$.

\end{lemma}

\begin{proof}
Assume $\vec{\phi}= \mathcal{H}_n\,\PP_n\vec{u}$.
We have already shown that $\vec{\phi}$ can be written in the form as given in (\ref{phihnm}).
Applying Lemma~\ref{rhorhotilde} and Lemma~\ref{lA} we find
\begin{eqnarray*}
\Phi_{h_{nm}(i)} &=&
\sum_{j=1}^{L(A_{\sigma_{n,m}(i)}0)} \big[ \mathcal{H}_n\,\rho_\Gnull{n}\big((m_j^{(h_{nm}(i))})^{-1}\big) \, \mathcal{H}_n^{-1} \, \vec{\phi} \big]_{h_n(\sigma_{m,n}(\Phi_{B_m}(i)))} \\
&& \qquad \qquad \qquad \qquad \qquad \qquad
\Big|_\beta m_j^{(h_{nm}(i))} \, \sigma_{R_i^{nm,1}}(B_m).
\end{eqnarray*}
Inserting $\vec{\phi}= \mathcal{H}_n\,\PP_n\vec{u}$ we find
\begin{eqnarray}
\nonumber
\Phi_{h_{nm}(i)} &=&
\sum_{j=1}^{L(A_{\sigma_{n,m}(i)}0)} \big[ \mathcal{H}_n\,\rho_\Gnull{n}\big((m_j^{(h_{nm}(i))})^{-1}\big) \, \PP_n\vec{u} \big]_{h_n(\sigma_{m,n}(\Phi_{B_m}(i)))} \\
\nonumber
&& \qquad \qquad \qquad \qquad \qquad \qquad
\Big|_\beta m_j^{(h_{nm}(i))} \, \sigma_{R_i^{nm,1}}(B_m)\\
\label{help}
&=&
\sum_{j=1}^{L(A_{\sigma_{n,m}(i)}0)} \big[ \rho_\Gnull{n}\big((m_j^{(h_{nm}(i))})^{-1}\big) \, \PP_n\vec{u} \big]_{\sigma_{m,n}(\Phi_{B_m}(i))} \\
\nonumber
&& \qquad \qquad \qquad \qquad \qquad \qquad
\Big|_\beta m_j^{(h_{nm}(i))} \, \sigma_{R_i^{nm,1}}(B_m).
\end{eqnarray}
Using again Lemma~\ref{lA} shows that $A_{\sigma_{n,m}(i)} = \sigma_{R_j^{nm,1}}(B_m)$.
Therefore the numbers $L(A_{\sigma_{n,m}(i)}0)$ in (\ref{help}) and $L(\sigma_{R_j^{nm,1}}(B_m)0)$ in (\ref{Pnm}) coincide.
Since also the matrices $m_j^{(h_{nm}(i))}$ in (\ref{help}) and $m_r^{(B_m,i)}$ in (\ref{Pnm}) coincide,  the identity $\Phi_{h_{nm}(i)} = (\PP_{nm} \vec{v})_i$ holds.
\end{proof}

\section{Hecke operators for $\Gnull{n}$}
\label{section3}

Proposition \ref{2 prop 1} allows us to introduce the following operators:

\begin{definition}
\label{2 def Hecke operators}
For $n \in \N$ define
\begin{equation}
\label{2 def hH}
\hH_{n,m}: S(n,\beta) \to S(n,\beta) \quad \mbox{by} \quad \hH_{n,m}(u)= \big(u \big|_0 B_m \big) \big|_0 \sum_j R_j^{nm,n}.
\end{equation}
\end{definition}

Obviously $\hH_{n,1}$ is the identity map $u \mapsto u$.

\begin{proposition}
\label{2 prop prim p}
For coprime $n,m \in \N$, the operator $\hH_{n,m}$ can be written as $\hH_{n,m} u = u \big|_0 \sum_{A\in X_m^\star} A$ with $X_m^{\star}$ given in (\ref{set Xn star}).
\end{proposition}

For reasons of simplicity we shall identify the operator $\hH_{n,p^{\nu}}$ with $\gcd(n,p)=1$ with the sum $\sum_{A\in X_{p^{\nu}}^{\star}} A$ when acting on Maass cusp forms.

To prove Proposition~\ref{2 prop prim p} we need the following

\begin{lemma}
\label{2 lemma p prim}
For $n,m \in \N$, $\gcd(n,m)=1$, there exists for every $j \in I_{nm,n}$   a matrix  $A(j)\in X_m^\star$ such that $B_m \, R_j^{nm,n} \, \big(A(j)\big)^{-1} \in \Gamma_0(n)$.
The map $j \mapsto A(j)$ is bijective.
\end{lemma}

\begin{proof}
Consider a system of representatives $\{R_j^{nm,n}\}_{j \in I_{nm,n}}$ as given in (\ref{2 system of representatives}).
Since $R_j^{nm,n} \in \Gnull{n} \subset \SL{\Z}$ there exists, as shown in \cite{Mu05}, an element $A(j) \in X_m^\star$ such that $\gamma_j = B_m {R_j}^{nm,n} \big(A(j)\big)^{-1} \in \SL{\Z}$.

We have to show that $\gamma_j \in \Gnull{n}$:
Since $R_j^{nm,n} \in \Gnull{n}$ we find
\[
\gamma_j = \Matrix {m}{0}{0}{1}\Matrix {a}{b}{nc}{d} {\Matrix {\alpha}{\beta}{0}{\delta}}^{-1}
=
\Matrix {m a \alpha ^{-1}}{-\beta a}{\frac{nc\delta}{m}}{\frac{-nc\beta + d\alpha}{m}} \quad \in \SL{\Z}.
\]
Since $\gcd(m,n)=1$ we conclude $m \mid c \delta$ and hence $\gamma_j \in \Gamma_0(n)$.

To show that the mapping $j \mapsto A(j)$ is injective we consider a $j^\prime \in I_{nm,n}$ such that $A(j^\prime) = A(j)$.
Put $g =R_j^{nm,n} \, \big( R_{j^\prime}^{nm,n} \big)^{-1}$.
By construction 
\[
\gamma_j \, \gamma_{j^\prime}^{-1}
=
B_m \, R_j^{nm,n} \, \big(A(j)\big)^{-1} \, \Big( B_m \, R_{j^\prime}^{nm,n} \, \big(A(j^\prime)\big)^{-1} \Big)^{-1}
=
B_m \, g \, B_m^{-1} \in \Gamma_0(n,m).
\]
Lemma~\ref{B7} shows that $g \in \Gnull{nm}$.
Hence $R_j^{nm,n}$ and $R_{j^\prime}^{nm,n}$ are representatives of the same right coset, implying $j=j^\prime$.

To finish the proof, we have to show surjectivity of the mapping $j \mapsto A(j)$. 
Since $\gcd(n,m)=1$ we know $[\SL{\Z}: \Gnull{n}] = [\Gnull{n}:\Gnull{mn}]$.
Hence the cardinalities of $I_n$, $I_{mn,n}$ and, since the map in (\ref{map Pn In}) is bijective, of $X_n^\star$ are equal.
\end{proof}

\begin{proof}[Proof of Proposition~\ref{2 prop prim p}]
By Definition~\ref{2 def Hecke operators} and Lemma~\ref{2 lemma p prim} we have
\[
\hH_{n,m}(u) = u \big|_0 \sum_j B_m\, R_j^{nm,n} = u \big|_0 \sum_j A(j) = u \big|_0 \sum_{A\in X_m^\star} A.
\]
\end{proof}

\begin{theorem}
\label{2 main theorem}
For fixed $n \in \N$ the family of operators $\hH_{n,m}$, $m \in \N$, satisfies
\begin{eqnarray*}
\hH_{n,p} \, \hH_{n,p^e}
&=&
\left\{\begin{array}{ll}
\hH_{n,p^{e+1}}
	& \mbox{for prime } p \mid n, e \in \N\\
\hH_{n,p^{e+1}} + p\,\Matrix{p}{0}{0}{p} \hH_{n,p^{e-1}} 
	& \mbox{for prime } p \nmid n, e>1,\\
\hH_{n,p^2} + (p+1)\,\Matrix{p}{0}{0}{p} \hH_{n,1} \qquad
	& \mbox{for prime } p \nmid n, e=1 \mbox{ and},
\end{array}\right.\\
\hH_{n,m} \, \hH_{n,m^\prime}
&=&
\hH_{n,mm^\prime}
\qquad\qquad\qquad\qquad\qquad\qquad\; \mbox{for } (m,m^\prime)=1.
\end{eqnarray*}
\end{theorem}

\begin{proof}
Take an $u \in S(n,\beta)$.
For $p$ prime we consider the three cases
\linebreak
$p \mid n$, $p \nmid n$ with $e=1$ and $p\nmid n$ with $e  \geq 2$ separately.

In the case $ p\mid n$ we find by using Definition~\ref{2 def Hecke operators} and Lemma~\ref{B4} 
\begin{eqnarray*}
\hH_{n,p^{e+1}} u
&=&
\big(u \big|_0 B_{p^{e+1}} \big) \big|_0 \sum_{l \in I_{np^{e+1},n}} R_l^{np^{e+1},n}\\
&=&
u \big|_0 \sum_{l \in I_{np^{e+1},n}} B_{p^{e+1}} \, R_l^{np^{e+1},n}\\
&=&
u \big|_0 \sum_{i \in I_{np,n}}\left( \sum_{j \in I_{np^e,n}} B_{p^e} \, R_j^{np^e,n} \right) \, B_p \, R_i^{np,n}\\
&=&
\hH_{n,p} \hH_{n,p^e} u.
\end{eqnarray*}

Consider next the case $p \nmid n$:
We know from Proposition~\ref{2 prop prim p} that the operators $H_{n,p^{e}}$ can be expressed for all $e \geq 1$ and for $p$ prime with $p \nmid n$ in terms of the elements $A \in X_{p^e}^\star$ as $H_{n,p^{e}} = \sum_{A \in X_{p^e}^\star} A$ for all $e \geq 1$.
Hence it is enough to determine the matrices in the set ${X_{p}}^\star \cdot {X_{p^e}}^\star = \{A \cdot B, A\in X_{p}^\star, B \in X_{p^e}^\star \} $.

For $e=1$ one finds
\begin{lemma}
\label{lemma e=1}
For $p$ prime 
\[
X_{p}^\star \cdot X_{p}^\star
=
X_{p^2}^\star
\cup \big\{\Matrix{p}{0}{0}{p} \big\}
\cup \big\{ \Matrix{p}{0}{0}{p} \, T^b; \, \,0\leq b\leq p-1 \big\}.
\]
\end{lemma}

\begin{proof}
Since for $p$ prime
\[
X_{p}^\star  = \{\Matrix{1}{b}{0}{p}; 0 \leq b \leq p-1 \} \cup \{\Matrix{p}{0}{0}{p}\}
\]
we find for $X_{p}^\star\cdot {X_{p}}^\star$:
\begin{eqnarray*}
\label{Eq p}
X_{p}^\star \cdot X_{p}^\star 
&=&
\big \{\Matrix{1}{b}{0}{p} \Matrix{1}{b^\prime}{0}{p};  0 \leq b \leq p-1,  0 \leq b^\prime \leq p-1 \big \}\\
&&
\cup \big\{\Matrix{1}{b}{0}{p} \Matrix{p}{0}{0}{1}; \, 0\leq b \leq p-1, \big\}
\cup \big\{ \Matrix{p^2}{0}{0}{1}\big\}\\
&&
\cup \big\{\Matrix{p}{0}{0}{1} \Matrix{1}{b^\prime}{0}{p}; \, 0\leq b^\prime \leq p-1 \big\}
\end{eqnarray*} 
and hence, splitting the sets into a part having coprime entries and noncoprime entries,
\begin{eqnarray*}
X_{p}^\star\cdot X_{p}^\star
&=&
\big\{\Matrix{1}{b^\prime + bp}{0}{p^2}; \quad 0\leq b \leq p-1,\quad 0\leq b^\prime \leq p-1 \big\}\\
&&
\cup \big\{ \Matrix{p}{b}{0}{p}: 1\leq b \leq p-1 \big\}
\cup \big\{ \Matrix{p^2}{0}{0}{1}\big\}\\
&&
\cup \big\{ \Matrix{p}{0}{0}{p} \big\}
\cup \big\{ \Matrix{p}{pb^\prime}{0}{p}; 0\leq b^\prime \leq p-1 \big\}.
\end{eqnarray*}
Using
\[
\big\{ \Matrix{p}{0}{0}{p} \big\} \cup \big\{ \Matrix{p}{pb^\prime}{0}{p}, 0\leq b^\prime \leq p-1 \big\}
=
\big\{ \Matrix{p}{0}{0}{p} \, \left( \mathbf{1} + \sum_{b=0}^{p-1} T^b \right) \big\}
\]
and, noticing that $b^\prime + bp$ above runs through all integers $\{0, \ldots, p^2-1\}$,
\begin{eqnarray*}
X_{p^2}^\star 
&=&
\big\{\Matrix{1}{b^\prime + bp}{0}{p^2}; \quad 0\leq b \leq p-1,\quad 0\leq b^\prime \leq p-1 \big\}\\
&&
\cup \big\{ \Matrix{p}{b}{0}{p}; 1\leq b \leq p-1 \big\} \cup \big\{ \Matrix{p^2}{0}{0}{1}\big\}
\end{eqnarray*} 
the proof of Lemma~\ref{lemma e=1} is finished.
\end{proof}

For $e \geq 2$ on the other hand one has

\begin{lemma}
\label{lemma e>2}
For $e \geq 2$ one has
\[
X_{p}^\star \cdot X_{p^e}^\star = X_{p^{e+1}}^\star \cup \bigcup_{l=0}^{p-1}\Matrix{p}{0}{0}{p} \Matrix{1}{l}{0}{1} X_{p^{e-1}}^\star.
\]
\end{lemma}

\begin{proof}
For
\begin{eqnarray*}
X_{p^e}^\star
&=&
\cup \big\{\Matrix{p^e}{0}{0}{1}\big\} \cup \big\{\Matrix{1}{b_0}{0}{p^e}; 0\leq b_0\leq p^e -1 \big\}\\
&&
\cup \bigcup _{j=1}^{e-1} \big\{\Matrix{p^j}{b_j}{0}{p^{e -j}}; 1\leq b_j\leq p^{e-j} -1, \gcd(b_j,p)=1) \big\}
\end{eqnarray*}
and
\[
X_p^\star = \big\{\Matrix{p}{0}{0}{1}\big\} \cup \big\{\Matrix{1}{a_0}{0}{p}; 0\leq a_0\leq p -1 \big\}
\]
we find for $X_{p}^\star \cdot X_{p^e}^\star$:
\begin{eqnarray*}
&& \!\!\!\!\!\!
X_{p}^\star \cdot X_{p^e}^\star \quad =\\
&&\quad
\big\{\Matrix{p^{e +1}}{0}{0}{1} \big\}
\cup \big\{ \Matrix{1}{b_0+a_0 p^e}{0}{p^{e +1}}; \, 0\leq a_0 \leq p -1, 0\leq b_0 \leq p^e -1 \big\}  \\
&&
\cup \bigcup_{j=1}^{e -1} \big\{ \Matrix{p^j}{b_j + a_0 p^{e -j}}{0}{p^{e +1-j}};\, 0\leq a_0\leq p-1, 1\leq b_j\leq p^{e-j} -1, \gcd(p,b_j)=1 \big\} \\
&&
\cup \bigcup_{j=1}^{e -1} \big\{ \Matrix{p^{j+1}}{pb_j}{0}{p^{e -j}}; \, 1\leq b_j\leq p^{e-j} -1, \gcd(p,b_j)=1 \big\}\\
&&
\cup \big\{ \Matrix{p^e}{a_0}{0}{p}; \, 0 \leq a_0 \leq p-1 \big\}
\cup \big\{\Matrix{p}{pb_0}{0}{p^e};\, 0\leq b_0\leq p^e-1 \big\}.
\end{eqnarray*}
Since $\tilde{b}_0 := b_0 + a_0 p^e$ takes the values $\{0, \ldots, p^{e+1}-1\}$ for $0 \leq b_0 \leq p^{e}-1$ and $0 \leq a_0 \leq p-1$ and $\tilde{b}_j := b_0 + a_0 p^{e -j}$ takes the values $\{1\leq l \leq p^{e-j+1}-1, \gcd(l,p)=1\}$ for $0 \leq a_0 \leq p-1$ and $1 \leq b_j \leq p^{e -j}-1$ we get
\begin{eqnarray*}
X_{p}^\star \cdot X_{p^e}^\star 
&=& 
X_{p^{e+1}}^\star
\cup \bigcup_{j=1}^{e -1} \big\{ \Matrix{p^{j+1}}{pb_j}{0}{p^{e -j}};\,  1\leq b_j\leq p^{e-j} -1, \gcd(p,b_j)=1 \big\} \\
&&
\cup \big\{ \Matrix{p}{pb_0}{0}{p^e};\, 0\leq b_0\leq p^e-1 \big\}\\
&=&
X_{p^{e+1}}^\star \cup \Matrix{p}{0}{0}{p}\,Y_e
\end{eqnarray*}
with
\begin{eqnarray*}
Y_e
&=&
\big\{ \Matrix{1}{b_0}{0}{p^e-1};\, 0 \leq b_0 \leq p^e -1 \big\}\\
&& 
\cup \bigcup_{j=1}^{e-1} \big\{ \Matrix{p^j}{b_j}{0}{p^{e-j-1}};\, 1 \leq b_j \leq p^{e-j}-1, (p,b_j)=1 \big\}.
\end{eqnarray*}
Let us consider the different subsets of $Y_e$ more in detail:
\begin{itemize}
\item
For $j=0$ we can write  $b_0$ satisfying $0 \leq b_0 \leq p^e -1$ uniquely as as $b_0=k_0+lp^{e-1}$ for some $0 \leq k_0 \leq p^{e-1}-1$ and some $0 \leq l \leq p-1$.
On the other hand we see that each $k_0+lp^{e-1}$ with $0 \leq k_0 \leq p^{e-1}-1$ and $0 \leq l \leq p-1$ satisfies $0 \leq k_0 + lp^{e-1} \leq p^e -1$.
Hence
\[
\big\{\Matrix{1}{b_0}{0}{p^{e-1}}; \, 0 \leq b_0 \leq p^e-1 \big\} =
\bigcup_{l=0}^{p-1} \big\{T^l\,\Matrix{1}{k_0}{0}{p^{e-1}}; \, 0 \leq k_0 \leq p^{e-1}-1 \big\}.
\]
\item
For $0 < j < e-1$ we can write $b_j$ satisfying $0 \leq b_j \leq p^{e-j} -1$ uniquely as $b_j=k_j+lp^{e-j-1}$ for some $0 \leq k_j \leq p^{e-j-1}-1$ and some $0 \leq l \leq p-1$.
On the other hand we see that each $k_j+lp^{e-j-1}$ with $0 \leq k_j \leq p^{e-j-1}-1$ and $0 \leq l \leq p-1$ satisfies $0 \leq k_j + lp^{e-j-1} \leq p^{e-j} -1$.
Furthermore the condition $\gcd(b_j,p)=1$ is equivalent to $\gcd(k_j,p)=1$.
Hence
\begin{eqnarray*}
&& \!\!\!\!
\bigcup_{j=1}^{e-1} \big\{ \Matrix{p^j}{b_j}{0}{p^{e-1-j}};\, 1\leq b_j\leq p^{e-j}-1, \gcd(p,b_j)=1 \big\}
\quad \\
&& =
\bigcup_{l=0}^{p-1} \big\{ T^l \, \Matrix{p^j}{k_j}{0}{p^{e-j-1}};\, 1 \leq k_j \leq p^{e-1-j}-1, \gcd(p,k_j)=1 \big\}.
\end{eqnarray*}
\item
For $j=e-1$ and $b_{e-1}$ with $0 \leq b_{e-1} \leq p-1$ finally we see $\Matrix{p^{e-1}}{b_{e-1}}{0}{1} = T^{b_{e-1}} \, \Matrix{p^{e-1}}{0}{0}{1}.$
Hence 
\[
\big\{\Matrix{p^{e-1}}{b_{e-1}}{0}{1}; \, 0 \leq b_{e-1} \leq p-1\big\}
=
\bigcup_{l=0}^{p-1} \big\{ T^l \, \Matrix{p^{e-1}}{0}{0}{1} \big\}.
\]
\end{itemize}
Summarizing the discussion above then shows that  $Y_e = \bigcup_{l=0}^{p-1} T^l \, X_{p^{e-1}}^\star$, which proves Lemma~\ref{lemma e>2}.
\end{proof}

Lemma~\ref{lemma e=1} and Lemma~\ref{lemma e>2} together with Proposition~\ref{2 prop prim p} prove the composition laws for the operators $\hH_{n,m}$ in Theorem~\ref{2 main theorem} for prime $p$.

To prove it for arbitrary $m,m'\in\mathbb{N}_0$ with $\gcd (m,m')=1$ we use Definition~\ref{2 def Hecke operators} and Lemma~\ref{B6} to get
\begin{eqnarray*}
\hH_{n,mm^\prime} u
&=&
\big(u \big|_0 B_{mm^\prime} \big) \big|_0 \sum_l R_l^{mm^\prime n,n}\\
&=&
u \big|_0 \sum_l B_{mm^\prime} \, R_l^{mm^\prime n,n}\\
&=&
u \big|_0 \sum_i \left( \sum_j B_{m^\prime} \, R_j^{m^\prime n,n} \right) \, B_m \, R_i^{mn,n}\\
&=&
\hH_{n,m} \hH_{n,m^\prime} u.
\end{eqnarray*}
This completes the proof of Theorem~\ref{2 main theorem}.
\end{proof}

From this Theorem we also get
\begin{corollary}
\label{cor1}
The operators $\hH_{n,m}$, $m \in \N$, commute.
\end{corollary}

\begin{proof}
We only have to show that $\hH_{n,p}\hH_{n,p^e}=\hH_{n,p^e}\hH_{n,p}$ for $p$ prime, $p\nmid n$ and arbitrary $e\in \mathbb{N}$.
This we show by induction on $e$.
For $e=1$ the operators trivially commute.
Assume they commute for all $1\leq e\leq N$.
Then one gets by Theorem~\ref{2 main theorem}
\begin{eqnarray*}
\hH_{n,p}\hH_{n,p^{N+1}}&=&\hH_{n,p}\big (\hH_{n,p}\hH_{n,p^N}-p \Matrix{p}{0}{0}{p}\hH_{n,p^{N-1}}\big)\\ 
&=&
\big (\hH_{n,p}\hH_{n,p^N}-p \Matrix{p}{0}{0}{p}\hH_{n,p^{N-1}}\big) \hH_{n,p}\\
&=&
\hH_{n,p^{N+1}}\hH_{n,p}.
\end{eqnarray*}
\end{proof}

\section{Proof of Theorem~\ref{main theorem}}
\label{section6}

The results of the preceeding sections lead to a simple interpretation of the operators $\tTT_{n,m}$ in (\ref{Hecke operator 1}) for arbitrary $n$ and $m$.

\begin{theorem}
\label{relation hH and tTT}
For $n,m \in \N$ and $u \in S(n,\beta)$ we have
\[
\tTT_{n,m} \mathcal{H}_n \PP_n \Pi_n u = \mathcal{H}_n \PP_n \Pi_n \hH_{n,m} u.
\]
\end{theorem}

\begin{figure}[h]
\label{diagram}
\setlength{\unitlength}{3ex}
\begin{picture}(12,7)
\put(1,1){$S(n,\beta)$}
\put(1,6){$S(n,\beta)$}

\put(2,5.5){\vector(0,-1){3.5}}
\put(0,3.5){$\hH_{n,m}$}

\put(4,1.2){\vector(1,0){4.5}}
\put(4.5,0.2){$\mathcal{H}_n\,P_n\,\Pi_n$}
\put(9,1){$\FE(n,\beta)$}

\put(4,6.2){\vector(1,0){4.5}}
\put(4.5,6.5){$\mathcal{H}_n\,P_n\,\Pi_n$}
\put(9,6){$\FE(n,\beta)$}

\put(10,5.5){\vector(0,-1){3.5}}
\put(10.5,3.5){$\tTT_{n,m}$}
\end{picture}

\caption{The relation between the operators $\hH_{n,m}$ in Definition~\ref{2 def Hecke operators} and $\tTT_{n,m}$ in Theorem~\ref{theorem 3}. Theorem~\ref{relation hH and tTT} shows that the diagram commutes.}
\end{figure}

\begin{proof}[Proof of Theorem~\ref{relation hH and tTT}]
For $u \in S(n,\beta)$ and $\hH_{n,m}$ the Hecke operators in Definition \ref{2 def Hecke operators} consider the Maass wave form $u^\prime = \hH_{n,m} u \in S(n,\beta)$.
Denote the corresponding vector valued Maass wave form by $\vec{u^\prime} = (u^\prime_i)_{i\in I_n} := \Pi_n \hH_{n,m} u$.
The vector valued period function $\PP_n\vec{u^\prime}$ then has the form
\begin{eqnarray*}
\big(\PP_n\vec{u^\prime}\big)_i (\zeta)
&=&
\int_0^\infty \eta(u^\prime_i,R_\zeta^\beta)
=
\int_0^\infty \eta\big((\Pi_n \hH_{n,m} u)_i,R_\zeta^\beta\big)\\
&=&
\int_0^\infty \eta\big((\hH_{n,m} u) \big|_0 R_i^{n,1} ,R_\zeta^\beta\big),
\qquad i\in I_n.
\end{eqnarray*}
Inserting the definition of the operator $\hH_{n,m}$ leads to
\begin{equation}
\label{P_n u'}
\big( \PP_n\vec{u^\prime} \big)_i (\zeta) = \int_0^\infty \eta\Big( \big( u|_0 B_m\sum_{j\in I_{nm,n}} R_j^{nm,n} \big) \big|_0 R_i^{n,1},R_\zeta^\beta \Big)
\end{equation}
and hence 
\[
\big( \PP_n\vec{u^\prime} \big)_i (\zeta) = \int_0^\infty \eta \Big( (v\big|_0 \sum_{j\in I_{nm,n}} R_j^{nm,n}) \big|_0 R_i^{n,1} , R_\zeta^\beta \Big)
\]
where $v=u\big|_0 B_m$ is an old Maass form for the group $\Gnull{nm}$ according to Proposition \ref{2 prop 1}.
Lemma~\ref{Lem} implies
\[
(\PP_n\vec{u^\prime})_i = \sum_{l\in\sigma_{m,n}^{-1}(i)} (\PP_{nm}\vec{v})_l.
\]
But by Lemma~\ref{Phi} we know that $\vec{\Phi} = \mathcal{H}_{nm} (\PP_{nm}\vec{v})$ if $\vec{\phi} = \mathcal{H}_n (\PP_n \Pi_n u)$ and hence we get
\[
(\PP_n\vec{u^\prime})_i
=
\sum_{l\in\sigma_{m,n}^{-1}(i)} (\PP_{nm}\vec{v})_l
=
\sum_{l\in\sigma_{m,n}^{-1}(i)} \Phi_{h_{nm}(l)}
\]
Lemma~\ref{h sigma} shows that $\sigma_{m,n}\big(h_{nm}(l)\big)=h_n(i)$ for all $l \in \sigma_{m,n}^{-1}(i)$.
Hence we write the equation above as
\[
(\PP_n\vec{u^\prime})_i
=
\sum_{l\in\sigma_{m,n}^{-1}(i)} \Phi_{h_{nm}(l)}
=
\sum_{l^\prime \in\sigma_{m,n}^{-1}(h_n(i))} \Phi_{l^\prime}
\]
Inserting the definition of the solution $\vec{\Phi}$ in Theorem \ref{theorem 2} we see that
\[
\sum_{l^\prime \in \sigma_{m,n}^{-1}(h_n(i))} \Phi_{l^\prime} = \left( \tTT_{n,m} \vec{\phi} \right)_{h_n(i)}.
\]
Since $\vec{\phi} = \mathcal{H}_n (\PP_n \Pi_n u)$ we finally get
\begin{equation}
\label{TH}
(\PP_n\vec{u^\prime})_i
=
\left( \tTT_{n,m} \mathcal{H}_n \PP_n \Pi_n u \right)_{h_n(i)}
\end{equation}
and hence 
\[
\mathcal{H}_n \PP_n\vec{u^\prime} = \tTT_{n,m} \mathcal{H}_n  \PP_n \Pi_n u.
\]
\end{proof}

\begin{proof}[Proof of Theorem~\ref{main theorem}]
Theorem~\ref{main theorem} follows now immediately from Theorem \ref{2 main theorem} and Theorem~\ref{relation hH and tTT} above.
\end{proof}

\begin{remark}
Proposition~3.13 in \cite{MM05}, establishing the relation between $\tTT_{n,m}$ and the Atkin-Lehner type ordinary Hecke operator $S(n,\beta) \to S(n,\beta)$; $u \mapsto u \big|_0 \sum_{A \in X_m} A$ for $m$ prime and $\gcd(n,m)=1$, is now a simple corollary of our Theorem~\ref{relation hH and tTT}.
\end{remark}

\begin{remark}
In his diploma thesis \cite{Fr05} M.\ Fraczek gives also an explicit form of the Fricke element when acting on vector valued period functions for the groups $\Gnull{n}$.
\end{remark}

\end{document}